\documentclass{amsart}
\usepackage[T1]{fontenc}

\usepackage{amssymb,amsfonts,amsmath}
\usepackage[all,arc]{xy}
\usepackage{enumerate}
\usepackage{mathrsfs}
\usepackage{graphicx}
\usepackage{hyperref}
\usepackage{todonotes}
\usepackage{xcolor}
\usepackage{pinlabel}
\usepackage{nicefrac}
\usepackage{bm}

\newtheorem{thm}{Theorem}[section]

\newtheorem{prop}[thm]{Proposition}
\newtheorem{lem}[thm]{Lemma}

\newtheorem{quest}[thm]{Question}

\newtheorem*{thm*}{Theorem}

\newtheorem{innercustomthm}{Theorem}
\newenvironment{customthm}[1]
  {\renewcommand\theinnercustomthm{#1}\innercustomthm}
  {\endinnercustomthm}

\theoremstyle{definition}
\newtheorem{defn}[thm]{Definition}
\newtheorem{defns}[thm]{Definitions}

\newtheorem{exmp}[thm]{Example}

\theoremstyle{remark}
\newtheorem{rem}[thm]{Remark}

\newcommand{\Z}{\mathbb{Z}}
\newcommand{\N}{\mathbb{N}}
\newcommand{\R}{\mathbb{R}}

\newcommand{\calA}{\mathcal{A}}
\newcommand{\calD}{\mathcal{D}}

\newcommand{\calC}{\mathcal{C}}
\newcommand{\calM}{\mathcal{M}}

\newcommand{\s}{\mathfrak{s}}
\renewcommand{\t}{\mathfrak{t}}
\renewcommand{\c}{\mathfrak{c}}

\DeclareMathOperator{\SO}{SO}
\DeclareMathOperator{\Fix}{Fix}
\DeclareMathOperator{\Branch}{Branch}
\DeclareMathOperator{\Arf}{Arf}
\DeclareMathOperator{\ind}{ind}

\DeclareMathOperator{\PD}{PD}

\title{Equivariant concordance of periodic 2-knots in $S^4$}
\author{Remy Bohm}
\address{Department of Mathematics, University of Texas, Austin, TX}
\email{remy.bohm@utexas.edu}

\begin{document}

\begin{abstract}
We show that the smooth equivariant concordance group of 2-knots in $S^4$ invariant under a linear $\Z/d\Z$ action is isomorphic to $\Z/2\Z$ for all $d \geq 2$.
This is in contrast to the non-equivariant case, in which all 2-knots are slice.
We construct a new invariant for these 2-knots, which we call periodic, and show that it fully classifies them up to equivariant concordance.
The invariant depends on a variation of the Arf invariant for null-homologous classical knots in arbitrary 3-manifolds with respect to a chosen spin structure.
Our proof also shows an identical classification for certain annuli in $S^1 \times B^3$ up to concordance rel. boundary.
\end{abstract}

\maketitle

\section*{Introduction}

Concordance of surfaces in various contexts has been an object of interest to topologists as early as the 1960s.
In 1965, Kervaire \cite{ker65} showed that all 2-spheres in $S^4$ are slice, meaning that they bound smoothly embedded 3-balls into $B^5$.
In 2015, Sunukjian \cite{sun15} generalized this to the statement that the concordance class of any oriented surface in a simply-connected 4-manifold is determined entirely by its genus and homology class.
Our work uses similar techniques as these earlier papers, and extends the theory of concordance of 2-spheres to a $\Z/d\Z$-equivariant setting.
In doing so, we end up proving a classification theorem for certain annuli in $S^1 \times B^3$ up to concordance which restricts to a product on the boundary.

Our goal is to examine \textit{equivariant concordance} of smooth 2-spheres that are symmetric with respect to an ambient order $d$ rotation $\rho$ acting on $S^4$ fixing an unknotted $S^2$, where the rotation restricted to the symmetric 2-sphere fixes exactly two points.
In analogy with the classical dimension, we call such 2-knots $d$-\textit{periodic}.
We say that these spheres are \textit{equivariantly slice} if they bound 3-balls into $B^5$ that are themselves equivariant with respect the extension of $\rho$ to $B^5$ that fixes an unknotted properly embedded 3-ball.
There is also a notion of periodic connected sum analogous to the connected sum of strongly invertible knots which turns the set of smooth equivariant concordance classes of periodic 2-knots into a group, denoted  $\calC^{d} _2$. 
This paper introduces a new invariant $\alpha_d$ of $d$-periodic 2-knots which completely determines their equivariant concordance class, simultaneously computing of all of these groups.

\begin{customthm}{\ref{thm: equi conc}}  
    For each positive integer $d$, there exists a $\Z/2\Z$-valued, additive concordance invariant $\alpha_d$ for $d$-periodic 2-knots $\tilde S$ such that:
    \begin{enumerate}
        \item if $\alpha_d(\tilde S) = 0$, then $\tilde S$ is $\Z/d\Z$-equivariantly slice,
        \item if $\alpha_d(\tilde S) = 1$, then $\tilde S$ is $\Z/d\Z$-equivariantly concordant to the $d$-twist spun trefoil, and
        \item every $d$-periodic 2-knot satisfies exactly one of the above two conditions.
    \end{enumerate}
    Thus the smooth equivariant concordance group $\calC^{d} _2$ is isomorphic to $\Z/2\Z$ for all $d$.
\end{customthm}

We prove our main theorem by first proving an identical classification for annuli in $S^1 \times B^3$ which are properly embedded such that both boundary components are copies of $S^1 \times \text{pt}$ inside of $\partial (S^1 \times B^3) = S^1 \times S^2$, which we call \textit{Montesinos annuli}.
We say that such an annulus is slice if it is concordant rel. boundary to a boundary parallel annulus in $S^1 \times B^3$.

\begin{customthm}{\ref{thm: annuli classification}}  
    There exists a $\Z/2\Z$-valued, additive relative concordance invariant $\alpha$ for Montesinos annuli $N$ in $S^1 \times B^3$ such that:
    \begin{enumerate}
        \item if $\alpha(N) = 0$, then $N$ is slice,
        \item if $\alpha(N) = 1$, then $N$ is concordant to the one-twist spun trefoil minus its spin axis, and
        \item every Montesinos annulus in $S^1 \times B^3$ satisfies exactly one of the above two conditions.
    \end{enumerate}
    Thus the smooth concordance group of Montesinos annuli $\calC _{2} ^\circ$ is isomorphic to $\Z/2\Z$.
\end{customthm}

The argument for these theorems go as follows.
Consider first a Montesinos annulus $N$.
We can associate to $N$ a ``corresponding'' 2-sphere $S$ in $S^4$ which intersect a ``removed'' 2-sphere $R$ in two points, such that removing a neighborhood of $R$ would turn $S$ into $N$ inside $S^1 \times B^3$.
We construct a Seifert 3-manifold for this corresponding sphere which intersect the removed sphere in a single arc, interpreting this as a sort of relative Seifert manifold for the annulus.
We then push this Seifert manifold into $B^5$, where it intersects a ``removed'' 3-ball $Q$ in a single arc, and then attempt to ambiently surger it to a 3-ball which intersect the removed 3-ball in a single unknotted arc.
This ambient surgery approach fails under certain conditions, which we show depends on an extension of the Arf invariant for knots to null-homologous knots in arbitrary 3-manifolds, detailed in Section \ref{sec: barf}.
Realizing this invariant as an annulus sliceness obstruction requires us to also build Seifert 4-manifolds for certain 3-manifolds in ambient 5-space with careful control of their intersection with an unknotted 3-sphere under specific Morse-theoretic conditions.

Once establishing this relative concordance group, we show that by lifting the ambient surgeries to the $d$-fold branched cover of $B^5$ over the removed 3-ball, we obtain a $\Z/d\Z$-equivariant slice ball for the branched cover of the sphere $S$ representing $N$ over its two intersection points with the removed sphere.
In this way, we interpret the removed sphere as the unknotted fixed set of a linear rotation on $S^4$, and the annulus concordance as the quotient of the periodic concordance by the action $\rho$.
This gives some intuition as to why the isomorphism type of these concordance groups, $\Z/2\Z$, does not depend on the period $d$ of the rotation.
That said, the following question remains open:

\begin{quest}
    Is the map from $\calC^{kd} _2$ to $\calC^{d} _{2}$ for $k \in \N$ induced by the inclusion of $kd$-periodic 2-knots into the set of $d$-periodic 2-knots is always an isomorphism?
    That is, does there exist a $kd$-periodic 2-knot which is not $\Z/kd\Z$-equivariantly slice, but is $\Z/d\Z$-equivariantly slice?
\end{quest}

This result fits into a body of work both on equivariant 2-knots and on surface concordance, in various categories.
Ruberman \cite{rub84} showed that 2-knots which are equivariant with respect to free involutions on $S^4$ are unique up to equivariant concordance.
Gabbard \cite{gab24} has also studied equivariant isotopy of 2-spheres that are also equivariant with respect to a rotation involution on $S^4$, but which are embedded in such a way that the rotation induces an orientation-reversing involution on the 2-sphere fixing a circle. 
Hellsten \cite{hel26} also has recent work on the classification of surfaces with boundary in simply-connected 4-manifolds up to concordance.

Related also is the work of McCooey \cite{mcc07} on topological concordance of $\Z/p\Z \times \Z/p\Z$ \textit{actions} on $S^4$ for $p$ prime, where he also found a $\Z/2\Z$-valued action concordance invariant.
McCooey's definition of concordance of actions implies concordance of their fixed sets, and so in proving his theorem, he shows that there are at most two topological concordance classes of pairs of 2-knots in $S^4$ intersecting in two points, often called \textit{Montesinos twins}, which are each the fixed set of a different commuting $\Z/p\Z$ action on $S^4$. His proof is surgery-theoretic, and the potential obstruction he identifies comes from the generator of the $\Z/2\Z$ part of $L_4(\Z[\Z \times \Z]) \cong \Z \oplus\Z/2\Z$.

Our work is substantially different from that of McCooey, reflecting the differing goals of his paper and ours, although we conjecture that our obstruction to smooth equivariant sliceness, restricted to the topological category and the correct class of periodic 2-knots, is equivalent to McCooey's.
Notably, McCooey's result only has implications for equivariant topological concordance within the class of 2-knots which are not only $\Z/p\Z$-equivariant, but are also themselves the fixed sets of $\Z/p\Z$ actions on $S^4$ for $p$ prime.
As mentioned before, his proof is fully surgery-theoretic, and it remains an open problem to find an example of a $\Z/p\Z \times \Z/p\Z$ action on $S^4$ for which his potential obstruction does not vanish.
In contrast, our theorem applies to all smooth periodic 2-knots of any period.
The proof is completely constructive, and we give explicit criteria for the vanishing of our obstruction and identify it as a direct extension of the Arf invariant of a knot, as well as examples where it does not vanish (although this does not immediately answer McCooey's question on the existence of nonlinear $\Z/p\Z \times \Z/p\Z$ actions on $S^4$).
Although the Arf invariant is typically a topological invariant, we work in the smooth category for the remainder of this paper.

There is also a family of literature on surface concordance in non-simply connected 4-manifolds, in which there are several known obstructions.
In 2018 Gabai gave a proof of the 4-dimensional lightbulb theorem \cite{gab20}, which gives conditions under which a sphere with a framed, geometrically dual sphere can be unknotted.
This result requires that the Freedman-Quinn invariant of the given sphere vanish in the case that the ambient 4-manifold has 2-torsion in its fundamental group.
The connection to the Freedman-Quinn invariant was formalized by Schneiderman and Teichner \cite{st22} after a counterexample to the lightbulb theorem in the case of ambient 2-torsion was constructed by Schwartz \cite{sch19}.
A similar result in the topological category had been proved by Freedman and Quinn \cite{fq90} showing that the Freedman-Quinn invariant obstructs concordance of 2-spheres, albeit with some errors which were corrected by Stong \cite{sto94}, and following Gabai's work, Miller \cite{mil21} gave a constructive proof.
In the case that a sphere in a non-simply connected 4-manifold has an unframed dual, Klug and Miller \cite{km21} extended work of Stong \cite{sto94} and Freedman and Quinn \cite{fq90} to show that the Kervaire-Milnor invariant of the sphere can be used to obstruct concordance.
Although we do not use these obstructions, we mention them in order to give a sense of scope of surface concordance invariants.

\subsection*{Outline}

In Section \ref{sec: prelim}, we establish some preliminaries of our setup for both Montesinos annuli in $S^1 \times B^3$ and periodic 2-knots.
In Section \ref{sec: barf}, we describe an extension of the classical Arf invariant to null-homologous knots in 3-manifolds, adapting work of Klug \cite{klu21}.
In Section \ref{sec: morse theory}, we outline some necessary results regarding embedded Morse theory with vertical boundary in the complement of a codimension two submanifold.
This requires the language of Morse theory for manifolds with boundary, as described by Borodzik-N\'emethi-Ranicki in \cite{bnr16}, as well as work of Borodzik-Powell in \cite{bp16}, and inspired by work of Schwartz in \cite{sch80} in the equivariant setting.
In Section \ref{sec: seif mfds}, we establish the existence of Seifert 3-manifolds for 2-knots while controlling the intersection of these manifolds with an unknotted 2-sphere intersecting our knot in two points. We then construct Seifert 4-manifolds for 3-manifolds in $S^5$ while controlling the intersection of these manifolds with an unknotted 3-sphere, given certain Morse theoretic conditions.
These both adapt a technique of Yanagawa in \cite{yana69}.
In Section \ref{sec: ann conc}, we use the extension in Section \ref{sec: barf} to define the $\alpha$ invariant for Montesinos annuli in $S^1 \times B^3$, and then give a procedure to ``slice'' these annuli via surgery on the Seifert 3-manifolds constructed in Section \ref{sec: seif mfds} under the assumption that $\alpha$ vanishes.
We also show that all Montesinos annuli for which $\alpha = 1$ are also mutually relatively concordant.
In Section \ref{sec: equi conc}, we use the results of the previous sections to define the invariant $\alpha_d$ for $d$-periodic 2-knots up to equivariant isotopy, and then lift all previous results to the equivariant setting.

\subsection*{Acknowledgments}
The author would like to thank Maggie Miller, Imogen Montague, Mark Powell, Anthony Conway, Evan Scott, Aru Mukherjea, and Seungwon Kim for their generous help and advice.
Special thanks to Malcolm Gabbard for suggesting this question, to Mark Powell for his helpful comments on an earlier draft of this paper, and to Imogen Montague, for pointing out and subsequently helping to fix a mistake in an earlier draft.

The author received partial support from the Simons Foundation collaboration grant on New Structures in Low-dimensional topology while writing this paper.

\section{Preliminaries on annuli and periodic 2-knots} \label{sec: prelim}

\subsection{Periodic 2-knots}

We begin with some background about smooth group actions on manifolds.
We specify to the case of $\Z/d\Z$ for this paper, but it is worth noting that there is a large body of work on 4-manifolds with various group actions; see, for example, Edmonds' survey \cite{edm18}.

\begin{defn}
    An action $\rho$ of $\Z/d\Z$ on the standard $m$-ball $B^m$ is \textit{linear} if it is conjugate to the action of an $m$-dimensional orthogonal representation of $\Z/d\Z$ on $B^m$.
    An action of $\Z/d\Z$ on $S^{m-1}$ is \textit{linear} if it is conjugate to the action of such a representation on $B^m$ restricted to the boundary.
\end{defn}

The actions we consider in this paper are the the linear rotation action $\sigma$ of $\Z/d\Z$ on $S^2$ fixing two points, and the linear rotation action $\rho$ of $\Z/d\Z$ on $B^5$ fixing an unknotted 3-ball such that its restriction to the boundary $S^4$ fixes an unknotted 2-sphere.
Restricting to linear actions simplifies the argument, although we hope to extend to nonlinear actions in future work.

Before proceeding, let us establish some notation for manifolds with cyclic group actions.
We will frequently need to pass between equivariant manifolds and their quotients by $\rho$.
Throughout this paper, we denote manifolds which carry an action by $\rho$ with tildes, to distinguish them from their quotients by $\rho$; e.g. $\tilde Y \subset \tilde B^5$ becomes $Y \subset B^5$ after taking the quotient by $\rho$.
(Because $\rho$ is always a linear cyclic action by rotation, the quotient of $\tilde B^5$ by $\rho$ is always diffeomorphic to $B^5$, and the quotient of $\tilde S^4$ by $\rho$ is always diffeomorphic to $S^4$.)
To recover the original manifold, we may take the $d$-fold cyclic branched cover over the image of the fixed set under the covering map.
For this reason, we refer to the image of $\Fix(\rho)$ under the covering map as the \textit{branch set} and denote it $\Branch(\rho)$.
Objects contained in the fixed set of $\rho$ will be denoted with hats, to distinguish them from their images under the quotient map (e.g. $\tilde Y \cap \Fix(\rho) = \hat A$, but in the quotient $Y \cap \Branch(\rho) = A$).

For convenience, we give a general definition for embeddings which are equivariant with respect to our desired type of action in all dimensions.

\begin{defn} \label{def: periodic}
    Let $\tilde M$ be a smooth $n$-manifold with a smooth action $\sigma$ of $\Z/d\Z$ on $\tilde M$, and let $\psi \colon \tilde M \to S^m$ be a smooth embedding for $m \geq n + 2$.
    We say that $\psi$ is \textit{$d$-periodic} if there exists an order $d$ linear rotation $\rho$ on $S^{m}$ fixing an unknotted smooth copy of $S^{m-2} \subset S^m$ such that $\psi(\tilde M)$ intersects $\Fix(\rho)$ transversely, and
    \[ \psi \circ \sigma = \rho \circ \psi. \]
    In the case that $\psi(\tilde M) \cap \Fix(\rho)$ is a single, smoothly embedded copy of $S^{n-2}$, or a single properly embedded copy of $D^{n-2}$ in the case that $\tilde M$ has boundary, we say that $\psi$ is \textit{simply $d$-periodic}.
\end{defn}

Some abuses of this notation:
where it is not explicitly relevant, we will often drop the prefix $d$- and describe embeddings just as \textit{periodic} or \textit{simply periodic}.
Furthermore, we will also often describe the \textit{image} of periodic embeddings as periodic; i.e., when we say that a 2-knot $\tilde S \subset S^4$ is periodic, we mean that there is an embedding $\psi \colon S^2 \to S^4$ whose image is $\tilde S$, and that $\psi$ is periodic for some period $d$.

Following this definition, we see that a 2-knot $\tilde S$ in $S^4$ is $d$-periodic if it satisfies
\[ \psi \circ \sigma = \rho \circ \psi,\]
with $\Fix(\rho)$ an unknotted 2-sphere.
In this case, since the restriction of $\rho$ to $\tilde S$ is orientation preserving, we see by the Lefschetz fixed point theorem that $S \cap \Fix(\rho)$ consists of exactly two points.
Such spheres arise naturally in common constructions of 2-knots: all $d$-twist spun knots can be constructed to have this symmetry about their spin axis.

\begin{rem}
    All periodic 2-knots are simply periodic as noted above, but the same is not true of surface knots in $S^4$.
    For instance, the unknotted torus in $S^4$ admits two different order 2 periodic symmetries: one in which one of the $S^1$ factors of the torus forms a meridian for $\Fix(\rho)$, and another that restricts to the hyperelliptic involution.
    The former action is free, while the latter fixes four points.
\end{rem}

Note that often, the intersection of a simply periodic manifold $\tilde M$ with the ambient fixed set can be knotted in $\tilde M$ (rel. boundary, if applicable).
This will be of particular importance to us.

\subsection{Montesinos annuli} \label{subsec: ann main defs}

As alluded to in the introduction and explained in detail in Section \ref{sec: ann conc}, we study the equivariant concordance of periodic 2-knots by studying the concordance of their quotients by $\rho$ relative to the fixed set of the action.
This concordance question is then equivalent to that of relative concordance of certain annuli in $S^1 \times B^3$.
Because of this, the majority of the paper takes place in the setting of these annuli, which we dub Montesinos annuli.

\begin{defn} \label{def: montesinos ann}
    A smooth \textit{Montesinos annulus} is a copy of $S^1 \times I$ smoothly and properly embedded in $S^1 \times B^3$ such that $\partial N \cong S^1 \sqcup S^1$ are mapped to disjoint copies of $S^1 \times \{\text{pt}\} \subset S^1 \times S^2 \cong \partial (S^1 \times B^3)$.
\end{defn}

We choose this name for these annuli in light of the following setup:
Consider $R$ a smoothly unknotted copy of $S^2$ inside of $S^4$, and let $S$ be a smooth 2-knot in $S^4$ intersecting $R$ transversely such that $S \cap R$ consists of exactly two points.
The pair $(R, S)$ is often called a \textit{Montesinos twin}, albeit one where one of the twins is unknotted.
Now consider removing a tubular neighborhood $\nu(R)$ from $S^4$, including neighborhoods of its intersection points with $S$.
What is left is an annulus embedded in $S^4 \setminus (D^2 \times S^2) \cong S^1 \times B^3$ exactly as described in Definition \ref{def: montesinos ann}.

With this in mind, we associate to each Montesinos annulus $N$ a \textit{corresponding sphere} $S$ and a \textit{removed sphere} $R$, such that $N \subset S^1 \times B^3$ is obtained by removing $\nu(R)$ from $S \subset S^4$.
It is important for this construction that we fix coordinates on the copy of $S^1 \times B^3$ in which the Montesinos annulus lives (and thus an identification of the boundary with $S^1 \times S^2$), so that the corresponding sphere is uniquely determined.
Otherwise, the Gluck map introduces an ambiguity when gluing back in the removed sphere.
Throughout the paper, we implicitly fix such an identification.

We say that a Montesinos annlulus is \textit{unknotted} if it is boundary parallel in $S^1 \times S^3$.
Note that the fundamental group of the complement of the unknotted annulus is $\Z^2$.

\begin{figure}[h]
    \labellist
    \small\hair 2pt
    \pinlabel {$\nu(R)$} at 46 300
    \pinlabel{$N$} at -10 100
    \pinlabel{$S$} at 273 100
    \pinlabel{$R$} at 368 300
    \endlabellist
    \centering
    \includegraphics[width=0.5\linewidth]{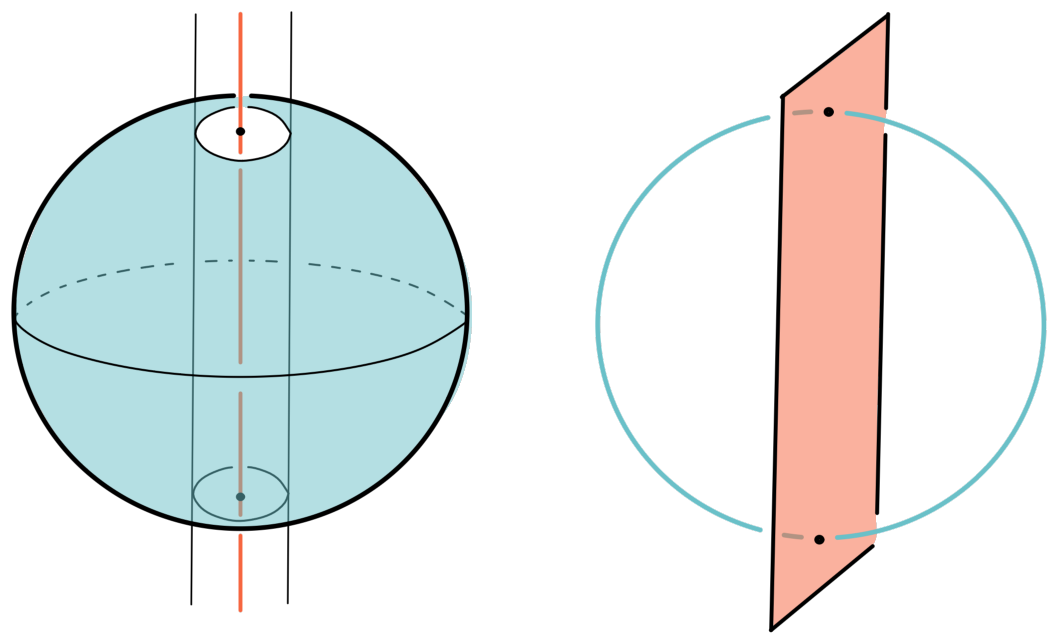}
    \caption{Two projections of an unknotted Montesinos annulus and its removed sphere.
    On the left, we see the annulus $N$ with its boundary components on $\partial (\nu(R))$.
    On the right and in a different 3-dimensional projection, we see $R$ locally as a plane, and one slice of the corresponding sphere $S$ for $N$ as a circle intersecting $R$ in two points.}
    \label{fig: ann vs sphere}
\end{figure}

\begin{figure}[h] 
    \labellist
    \small\hair 2pt
    \pinlabel{$T$} at 70 40
    \pinlabel{$A$} at 18 90
    \endlabellist
    \centering
    \includegraphics[width=0.8\linewidth]{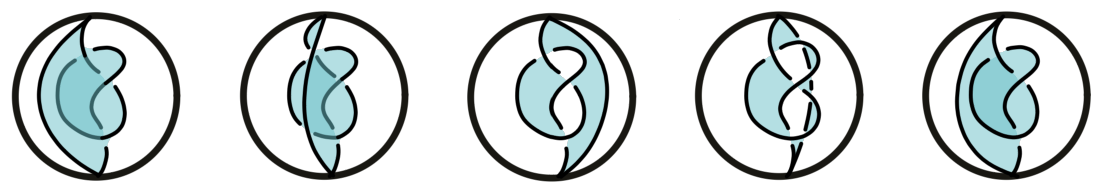}
    \caption{A movie of the one-twist spun trefoil, illustrated by twisting an the spin axis $A$ on the boundary while fixing the trefoil tangle $T$. The Seifert surfaces shown patch together to form a 3-ball in which the complement of $A$ is fibered by punctured tori with the trefoil monodromy.}
    \label{fig: spinning trefoil}
\end{figure}

\begin{exmp} \label{ex: N_T}
    Let $S_T$ be the one-twist spun trefoil 2-knot inside $S^4$ (not up to isotopy), and let $R$ be an unknotted copy of $S^2$ representing its spin axis.
    Remove a neighborhood $\nu(R)$ from $S^4$, including neighborhoods of its intersection with $S_T$.
    This leaves us with a Montesinos annulus inside of $S^1 \times B^3$, which we will call $N_T$.

    The key point to illustrate is that although the sphere $S_T$ is isotopic to the unknot in $S^4$ (see Zeeman \cite{zee65}), the annulus $N_T$ is not isotopic to the unknotted annulus rel. boundary.
    It is distinguished by the fundamental group of its complement in $S^1 \times B^3$, which is isomorphic to the direct product of the trefoil knot group with $\Z$.

    Note also the following construction, illustrated in Figure \ref{fig: spinning trefoil}, which will be important later.
    Following Zeeman, we can construct a 3-ball $B$ in $S^4$ bounded by $S_T$ by twist-spinning a relative Seifert surface $F$ for the trefoil tangle $T$.
    Close inspection of this construction shows that the isotopy type of the spin axis $A$ inside of this 3-ball is itself also a trefoil.
    The Seifert surfaces for $T$ together form a genus one fibration of the complement of $A$ in $B$ which glue together with precisely the trefoil monodromy, as shown in Figure \ref{fig: spinning trefoil}.
    Hence the intersection of $B$ with the removed sphere $R$ for $N_T$ is also a trefoil tangle.
\end{exmp}

\begin{defn} \label{def: ann conc}
    Two Montesinos annuli $N$ and $N'$ are smoothly \textit{concordant} if there exists $E$ a copy of $S^1 \times I \times I$ smoothly and properly embedded into $S^1 \times B^3 \times I$ such that $E|_{S^1 \times I \times \{0\}} = N$, $E|_{S^1 \times I \times \{1\}} = N'$, and $E|_{S^1 \times \{0,1\} \times I} = S^1 \times \{p_+, p_-\} \times I$ for some fixed points $p_+, p_-$ in $S^1 \times S^2$. 

    Equivalently, $N$ and $N'$ are concordant if their corresponding spheres $S$ and $S'$ are concordant in $S^4 \times I$ such that the intersection of the concordance with $R \times I \subset S^4 \times I$ consists of two arcs $\{p_+\} \times I, \{p_-\} \times I$ which are embedded in $R \times I$ by a level preserving map.
    We say that an annulus $N$ is \textit{slice} if it is concordant to the unknotted annulus. 
\end{defn}

\begin{rem} \label{rem: ann slice}
    Note that the kind of relative concordance for annuli we consider requires that both boundary components of the annulus stay fixed throughout the whole concordance, and are not allowed to vary in their isotopy class.
    Because of this restriction, we may obtain a different characterization of slice annuli adapted to the setting of $B^5$.
    Parametrize a neighborhood of $B^5$ near its boundary with coordinates $(x,y,z,w,t)$ such that $t \geq 0$.
    Consider $S^1 \times B^3 \times I$ sitting inside of $B^5$ such that $S^1 \times B^3 \times \{0\}$ is contained in $\partial B^5$, and so that projection onto the $I$ factor agrees with projection onto the $t$ coordinate.
    If $N$ is slice, then we may consider it as contained in $S^1 \times B^3 \times \{0\}$, and the unknotted annulus $U$ as contained in $S^1 \times B^3 \times \{1\}$.
    By filling in a 3-ball bounded by the corresponding sphere for $U$ lying in the $t=1$ slice (which necessarily intersects the removed sphere $R$ in that slice in an unknotted arc) and perturbing the 3-ball slightly (in a way formalized in Lemma \ref{lem: pushing into b5}), we see that the corresponding sphere $S$ for $N$ bounds a 3-ball into $B^5$ whose intersection with an unknotted copy $Q$ of $B^3$ has exactly one critical point with respect to projection onto the $t$ coordinate.
    We call such a slice ball minus a neighborhood of $Q$ a \textit{slice cylinder} for $N$.
    This characterization is the one we will work with most often, and many of our constructions and proofs take place in this parametrization.
\end{rem}

Throughout this paper, given a Montesinos annulus $N$, $R$ will always stand for its removed sphere, and $S$ for its corresponding sphere.
For brevity, we will also often drop the word ``Montesinos'' from Montesinos annulus, where it is clear from context.

\subsection{Concordance group structure}

\subsubsection{Concordance of annuli}

To define a group structure on the set of Montesinos annuli up to concordance, we require a notion of connected sum, which in turn necessitates some extra data.
We take inspiration from the connected sum for strongly invertible classical knots first defined by Sakuma in \cite{sak86} in our definition.
Throughout this subsection, $N$ is a Montesinos annulus, and $S$ and $R$ are its corresponding and removed spheres, respectively.

\begin{defn}
    A \textit{direction} for an annulus $N$ is a choice of orientation on $R$ together with a choice of labeling of the two boundary components of $N$ with the labels $c_+$ and $c_-$.
\end{defn}

We call an annulus with a choice of direction a \textit{directed} annulus.
The concordance of annuli given in Definition \ref{def: ann conc} can then be refined to a notion of \textit{directed concordance}, in which we require that the directions of the two annuli on either end of the concordance agree.

\begin{defn} \label{def: stacking}
    Let $N$ and $N'$ be two directed annuli, and let $R$ and $R'$ denote the removed spheres in their respective copies of $S^4$.
    Then we define their \textit{stack}, denoted $N \natural N'$, to be the annulus obtained as follows:
    From the copy of $S^4$ containing $N$, remove a ball such that $c_+$ lies in the new $S^3$ boundary component. 
    From the copy of $S^4$ containing $N'$, remove a ball such that $c_-'$ lies in the new $S^3$ boundary component.
    Then glue these two copies of punctured $S^4$ together along their boundaries such that:
    \begin{enumerate}
        \item the intersections of $R$ and $R'$ with the deleted neighborhoods are glued together along their common boundary with opposite orientations, so that the resulting sphere $R \# R'$ has a consistent orientation, and
        \item $N$ and $N'$ are glued so that $c_+$ is identified with $c_-'$ with opposite orientations.
    \end{enumerate}
    The result is an annulus $N \natural N'$ whose corresponding sphere is the connected sum $S \# S'$ which intersects the connected sum of removed spheres $R \# R'$ in two points.
\end{defn}

\begin{figure}[h]
    \labellist
    \small\hair 2pt
    \pinlabel{$R \# R'$} at 155 -10
    \pinlabel{$N$} at 130 50
    \pinlabel{$N'$} at 110 255
    \pinlabel{$S^3$} at 350 150
    \pinlabel{$c_+$} at 230 150
    \pinlabel{$c_-'$} at 140 180
    \endlabellist
    \centering
    \includegraphics[width=0.5\linewidth]{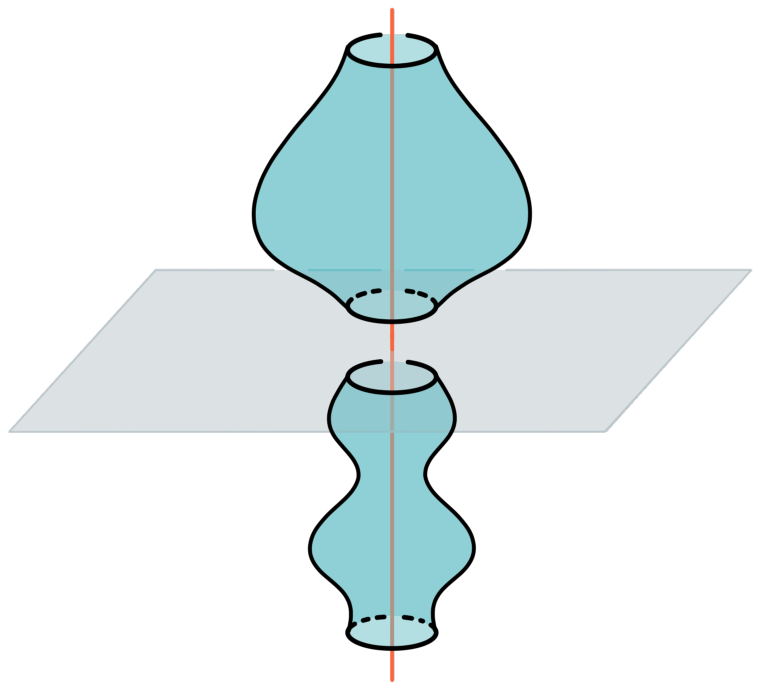}
    \caption{A schematic of one annulus stacked on top of another.}
    \label{fig: stacking}
\end{figure}

As in the case with strongly invertible classical knots, annulus connected sum is not usually commutative; that is, $N \# N'$ may not be isotopic to $N' \# N$.
In comparison with the classical dimension, note Di Prisa \cite{dip23} has shown that the equivariant concordance group of strongly invertible knots is not itself commutative, and so equivariant connected sum is not even commutative up to concordance. 
However, we will see that this is not the case for the concordance group for Montesinos annuli $\calC_2 ^\circ$.

\begin{defn}
    Let $C_2 ^\circ$ be the set of directed Montesinos annuli.
    The \textit{relative concordance group} $\calC_2 ^\circ $ of Montesinos annuli is the quotient of the set $C_2 ^\circ$ by directed concordance, with group multiplication given by stacking.
\end{defn}

\subsubsection{Equivariant concordance}

The situation for periodic 2-knots is almost identical.
We reuse much of the same language.

\begin{defn}
    A \textit{direction} for a $d$-periodic 2-knot $(\tilde S, \rho)$ is a choice of orientation on $\Fix(\rho)$ together with a choice of labeling of the two points in $\Fix(\rho) \cap \tilde S$ with the labels $\hat p_+$ and $\hat p_-$.
\end{defn}

We call a $d$-periodic 2-knot with a choice of direction a \textit{directed} $d$-periodic 2-knot.

\begin{defn} \label{def: equi conn sum}
    Let $(\tilde S, \rho)$ and $(\tilde S', \rho')$ be two directed $d$-periodic 2-knots.
    Then we define their \textit{equivariant connected sum}, denoted $(\tilde S\# \tilde S', \rho\#\rho')$ or just $\tilde S \#_{\Z/d\Z}\tilde S'$, to be the $d$-periodic 2-knot obtained as follows:
    First, remove an equivariant neighborhood of $p_+$ in $S^4$ from the copy of $S^4$ containing $\tilde S$, and an equivariant neighborhood of $p_-'$ from the copy of $S^4$ containing $\tilde S'$.
    Then glue these two copies of $S^4$ together along their boundaries such that:
    \begin{enumerate}
        \item the intersections of $\Fix(\rho)$ and $\Fix(\rho')$ with the deleted neighborhoods are glued together along their common boundary with opposite orientations (so that the union of $\Fix(\rho)$ with $\Fix(\rho')$ has a consistent orientation given by the chosen orientations on each of its components),
        \item $\rho$ and $\rho'$ agree on the common boundary of the punctured 4-spheres, and
        \item the punctured copies of $\tilde S$ and $\tilde S'$ are glued together along their boundaries.
    \end{enumerate}
    The resulting 2-knot $\tilde S\#_{\Z/d\Z} \tilde S'$ will then be $d$-periodic under the action given by the actions $\rho$ and $\rho'$ on their constituent components.
\end{defn}

For the following definition, extend the linear order $d$ rotation $\sigma$ on $S^2$ to $S^2 \times I$ to the action $\tau$ such that $\tau(x,t) = (\sigma(x), t)$.

\begin{defn}
    An \textit{equivariant concordance} of $d$-periodic 2-knots is an embedding $\Psi \colon S^2 \times I$ into $S^4 \times I \subset B^5$ such that:
    \[ \rho \circ \Psi = \Psi \circ \tau \]
    satisfying the condition that 
    \[\Psi|_{\Fix(\tau)}(x,t) = (\psi(x), t)\] 
    for some map $\psi \colon \Fix(\sigma) \to \Fix (\rho)$.
\end{defn}

\begin{rem} \label{rem: equi conc def}
    This second condition amounts to requiring that our equivariant concordances be level preserving on the fixed sets of the action.
    This is equivalent to the boundary conditions in Definition \ref{def: montesinos ann} for concordances of Montesinos annuli, and is necessary in order for the equivariant connected sum operation defined above to act as a group multiplication for equivariant concordance classes.
    In principle, one could study equivariant concordances or equivariant sliceness without this condition, but it is not clear that the set of equivalence classes under a looser definition of equivariant concordance actually forms a group.
    Also, removing this requirement complicates many Morse-theoretic arguments which are central to the results in this paper, and it is not immediately obvious how to recover these results without them.
\end{rem}

We will work primarily with $\tilde B^5$ and think of the equivariant concordance as taking place inside an embedded copy of $\tilde S^4 \times I$ such that the actions on each $\tilde S^4$ are the restriction of a global order $d$ rotation on $\tilde B^5$ fixing a properly embedded unknotted 3-ball.
If $(\tilde S, \rho)$ and $(\tilde S', \rho')$ are directed, we say that a concordance between them is directed if the orientations on $\Fix(\rho)$ and $\Fix(\rho')$ agree.
Directed equivariant concordance is then an equivalence relation on the set of directed periodic 2-knots.

\begin{defn}
    Let $C_2 ^d$ be the set of directed $d$-periodic 2-knots.
    The \textit{equivariant concordance group} $\calC_2 ^d$ of $d$-periodic 2-knots is the quotient of the set $C_2 ^d$ by directed equivariant concordance, with group multiplication given by equivariant connected sum.
\end{defn}

It is worth mentioning at this point that the condition of equivariant concordance for 2-knots is well-defined up to equivariant isotopy.
For, if $\tilde S$ is equivariantly concordant to $\tilde S'$ and $\tilde S''$ is equivariantly isotopic to $\tilde S$, then we may find an equivariant concordance between $\tilde S'$ and $\tilde S''$ by tracing the equivariant isotopy from $\tilde S''$ to $\tilde S$ and then finishing with the concordance from $\tilde S$ to $\tilde S'$.

\section{An extension of the Arf invariant} \label{sec: barf}

In this section, we build a necessary invariant that we will adapt later to the setting of annulus and periodic concordance.

\subsection{Background}

The Arf invariant of a quadratic form was introduced by Arf \cite{arf41} as a purely algebraic tool, and was identified as an invariant that could be associated to a classical knot by Robertello \cite{rob65}.
His paper describes three invariants of knots in $S^3$, then proves they are all equivalent.
For our purposes, it will be helpful to recall all three.

Before stating them, we remind the reader of the following definition.
An $n-2$-manifold $N$, possibly with boundary, properly embedded in an $n$-manifold $M$ is called \textit{characteristic} if $\PD([N]) = w_2(M) \in H^2(M;\Z/2\Z)$.
In the case of a surface $F$ in a 4-manifold $X$, there is an equivalent condition:
$F$ is characteristic in $X$ if and only if for all $x \in H_2(X;\Z/2\Z)$, we have that $[F] \cdot x \equiv x \cdot x$ mod 2, where $\cdot$ denotes algebraic intersection number.

A note on notation: in this section, $\sigma(X)$ refers to the signature of a 4-manifold $X$, and $\mu(Y, \s)$ refers to the Rokhlin invariant of a 3-manifold $Y$ with respect to the spin structure $\s$.
Where the 3-manifold has a unique spin structure, as is the case for an integer homology sphere $\Sigma$, for example, we use just the notation $\mu(\Sigma)$.

\begin{defns}[Robertello]\label{def: robertello}
    Let $K$ be an oriented knot in $S^3$.
    \begin{enumerate}
        \item Let $F$ be a Seifert surface for $K$ of genus $g$.
        Choose a symplectic basis for $H_1(F;\Z/2\Z)$ and let $I_F\colon H_1(F;\Z/2\Z) \times H_1(F;\Z/2\Z) \to \Z/2\Z$ be the associated Seifert form, considered as a $2g\times2g$ matrix over $\Z/2\Z$ in the chosen basis.
        Let $q_F \colon H_1 (F;\Z/2\Z) \to \Z/2\Z$ be given by:
        \[ q_F (v) = I_F(v,v) \]
        Then $q_F$ is a quadratic form.
        \textit{Arf invariant} of $K$ is the Arf invariant of $q_F$.
        \item Let $X$ be a smooth 4-manifold with boundary $S^3$ such that $K$ bounds a smooth, properly embedded characteristic disk $D$ into $X$.
        Define the quantity $\varphi(K) \in \Z/2\Z$ by:
        \[ \varphi(K) := \frac{D \cdot D - \sigma(X)}{8} \mod 2 \]
        where $D \cdot D$ here represents the normal Euler number of $D$ relative to the Seifert framing on $K$.
        \item We say $K$ and $J$ are \textit{band pass equivalent} if there exists a sequence of planar isotopies, Reidemeister moves, and band passes (see Figure \ref{fig: band pass}) taking a diagram of $K$ to a diagram of $J$.
        Band pass equivalence is an equivalence relation on the set of knots in $S^3$, and divides knots into two classes: those band pass equivalent to the unknot, and those band pass equivalent to the trefoil.
        We see therefore that the \textit{band pass equivalence class} is an invariant of $K$.
    \end{enumerate}
\end{defns}

\begin{thm}[Robertello] \label{thm: arf 0 band pass}
    The three definitions given above are all equivalent to each other.
    Specifically, $K$ has Arf invariant 0 if and only if $\varphi(K)$ vanishes, and equivalently, if and only if it is band pass equivalent to the unknot.
\end{thm}

Most treatments of the Arf invariant of knots use the definition (1); for instance, this is the formal definition given in \cite{sav12}, Chapter 9.
However, there are also several extensions and related invariants in 3- and 4-manifold topology that build on the second definition.
An excellent overview of the progression of extensions of the Arf invariant can be found in Klug's thesis \cite{klu21}.

For the extension we wish to define, we build on work of Klug, and in fact recover a special case of an extension of the Brown invariant defined in his paper.
The model equation for Arf-like invariants in the smooth category is as follows: given a characteristic surface $F$ smoothly and properly embedded in a smooth 4-manifold $X$ with $\partial F = K$ and $\partial X$ an integer homology 3-sphere $\Sigma$ with Rokhlin invariant $\mu(\Sigma)$, we have:
 \[ \Arf(F) + \Arf(K) \equiv \frac{F\cdot F - \sigma(X)}{8} + \mu(\Sigma) \mod 2. \]
Here $\Arf(F)$ denotes the Arf invariant of the Rokhlin form on $F$ (this term vanishes if $F$ is a disk, and so we won't address it further), and $F \cdot F$ denotes the normal Euler number of $F$ relative to the Seifert framing on $K$.

The goal of our invariant is to remove the requirement that the 3-manifold in which $K$ lives be an integer homology sphere, at the cost of requiring extra information.

\subsection{The Arf invariant of a null-homologous knot in a spin 3-manifold} \label{subsec: barf}

Following Kirby-Taylor \cite{KT90}, Klug \cite{klu21} showed that to extend the Arf invariant to the case of null-homologous knots in 3-manifolds, we need the data of a spin characterization.
Our contribution is to note that there exists a canonical spin characterization in the case that the knot is null-homologous over $\Z$, which allows for a more natural invariant with a clear concordance interpretation.

\begin{defn}
    Let $M$ be an $n$-manifold possibly with boundary, and let $N$ be a properly embedded characteristic codimension 2 submanifold.
    A \textit{spin characterization} $\c$ of $N$ in $M$ is a spin structure on $M \setminus \nu(N)$ that does not extend over $\nu(N)$.
\end{defn}

Just as in the case of spin structures, $H^1(M; \Z/2\Z)$ acts freely and transitively on the set of spin characterizations of $N$ in $M$.

We will be interested in the case where the characteristic submanifold is a knot $K$ in a 3-manifold $Y$ that is null-homologous over $\Z$.
This is strictly stronger than being dual to $w_2(Y)$: since all 3-manifolds are spin, $w_2(Y) = 0$, and so it suffices for $K$ to be null-homologous over $\Z/2\Z$.
However, the independence of the Arf invariant from a choice of spin characterization does not hold in the case that $K$ represents the zero class in homology with $\Z/2\Z$ coefficients but a nontrivial class in homology with $\Z$ coefficients.
From here on, when we say $K$ is null-homologous, we mean that $[K] = 0$ in $H_1(Y;\Z)$. 

\subsubsection{The induced characterization from a spin structure} \label{subsubsec: induced char}
Given a spin structure $\s$ on $Y$, when $K$ is null-homologous in $Y$ there is a unique induced spin characterization of $K$ in $Y$ that can be obtained as follows.
Consider a characterization $\c$ of $K$, and restrict $\s$ to the complement of a neighborhood of $K$.
Then since $\s$ and $\c$ are both spin structures on $Y \setminus \nu(K)$, their difference $[\s - \c]$ is an absolute class in $H^1(Y \setminus \nu(K); \Z/2\Z)$.
Because $\Z/2\Z$ is a field, $H^1(Y \setminus \nu(K); \Z/2\Z) \cong H_1(Y \setminus \nu(K); \Z/2\Z)$ by the universal coefficient theorem, and so we may interpret $[\s - \c]$ as a 1-cycle over $\Z/2\Z$ in $Y$.
Note that because $\s$ and $\c$ disagree on $\partial \nu(K)$ and in particular on a meridian of $K$, this cycle is never trivial in homology.
We can therefore make the following definition:
\begin{defn} \label{def: induced spin char}
    The spin characterization $\c_\s$ of $K$ in $Y$ \textit{induced} by $\s$ is the unique spin characterization $\c$ such that $[\s-\c_\s]$ is the image of the meridian of $K$ in $H^1(Y \setminus \nu (K); \Z/2\Z)$ under the map given by the universal coefficient theorem.
\end{defn}

From this we also get a preferred Seifert surface for $K$ via Poincar\'e-Lefschetz duality.
The group $H^1(Y \setminus \nu(K); \Z/2\Z)$ dual to $H_2(Y/\nu(K), \partial (\nu(K)); \Z/2\Z)$, which is isomorphic to $H_2(Y, K; \Z/2\Z)$ by excision.
Let $\beta$ be the image of $[\s - \c_\s]$ under this composition of isomorphisms.
Then $\beta$ lifts to an integral class $\delta$ because it is dual under the intersection pairing on $Y$ to the meridian of $K$, which is non-torsion.
Hence $\delta$ is representable by a Seifert surface $F_\s$ for $K$ which is determined up to its integer homology class relative to $K$.

\begin{rem}
    For any characteristic submanifold of a spin manifold given a spin structure, we have an induced spin characterization of the submanifold given analogously.
    Thus Definition \ref{def: induced spin char} makes sense in all dimensions.
    Note also, as mentioned earlier, that the set of spin characterizations is in bijection with the first cohomology group of the ambient manifold with $\Z/2\Z$ coefficients.
    Klug \cite{klu21} explores extensions of the Arf invariant compatible with other characterizations.
    However, for our purposes, we will only need this canonical one.
\end{rem}

\subsubsection{Definition of the invariant}

For brevity, we make the following auxiliary definition.

\begin{defn} \label{def: char filling}
    Let $(Y, \s)$ be a spin 3-manifold and let $K$ be a null-homologous knot in $Y$.
    A smooth \textit{characteristic filling} of $(K, Y, \s)$ is a pair $(X, D)$ where $X$ is a smooth 4-manifold with $\partial X = Y$ and $D$ is a characteristic disk in $X$ with $\partial D = K$ such that there exists a spin characterization of $D$ in $X$ extending the canonical spin characterization $\c_\s$ of $K$ in $Y$.
\end{defn}

We can now define our extension of the Arf invariant.

\begin{defn} \label{def: barf}
    Let $(Y, \s)$ be a spin 3-manifold and let $K$ be a null-homologous knot in $Y$.
    Let $(X, D)$ be a characteristic filling of $(K, Y, \s)$.
    We define the \textit{Arf invariant} of $K$ in $Y$ with respect to the spin structure $\s$, denoted $\Arf(K,Y,\s)$, to be:
    \[ \Arf(K,Y,\s) := \frac{D \cdot D - \sigma(X)}{8} + \mu(Y, \s) \mod 2 \ \in \Z/2\Z.\]
\end{defn}

It is not hard to see that this Arf invariant is well-defined:

\begin{prop} \label{prop: barf kys well-def}
    $\Arf(K,Y,\s)$ does not depend on the choice of smooth characteristic filling.
\end{prop}
\begin{proof}
    Let $(X_1, D_1)$ and $(X_2, D_2)$ be two smooth characteristic fillings of $Y$, and let  $\c_i$ be spin characterizations of $D_i$ in $X_i$ restricting to $\c_\s$ on the boundary. 
    We may then glue $X_1$ to $-X_2$ along their boundary so that $D_1 \cup _K -D_2$ forms a sphere, which is characteristic in $X_1 \cup_Y -X_2$ precisely because the spin characterizations on each $D_i$ glue together along $\c$ to give a spin characterization of $D_1 \cup _K -D_2$.
    The existence of such a spin characterization implies that $D_1 \cup _K -D_2$ is characteristic because the second Stiefel-Whitney class of the complement vanishes exactly when $D_1 \cup _K -D_2$ is dual to that class.
    Thus we may apply a theorem of Kervaire-Milnor \cite{KM61} to conclude that:
    \[ (D_1 \cup _K -D_2) \cdot (D_1 \cup _K -D_2) \equiv \sigma(X_1 \cup_Y -X_2) \mod 16.\]
    By Novikov additivity \cite{nov70}, $\sigma(X_1 \cup_Y -X_2) = \sigma(X_1) - \sigma(X_2)$, and since $D_i \cdot D_i$ was computed with respect to the Seifert framing on $K$, we have that:
    \[ (D_1 \cup _K -D_2) \cdot (D_1 \cup _K -D_2) = D_1 \cdot D_1 - D_2 \cdot D_2.\]
    Thus the difference of the Arf invariants computed using $(X_1, D_1)$ and $(X_2, D_2)$, after removing the Rokhlin invariant term that simply cancels with itself, is:
    \begin{align*}
        &\frac{D_1 \cdot D_1 - \sigma(X_1)}{8} - \frac{D_2 \cdot D_2 - \sigma(X_2)}{8} \\
        &= \frac{(D_1 \cup _K -D_2)\cdot (D_1 \cup _K -D_2) - (\sigma(X_1) - \sigma(X_2))}{8} \equiv 0 \mod 2.
    \end{align*}
    Because $\Arf(K, Y, \s)$ is valued in $\Z/2\Z$, the result follows.
\end{proof}

\subsection{Spin characteristic fillings}

In this subsection, we give an algorithm to produce smooth characteristic fillings as described in Definition \ref{def: char filling}.
Specifically, the algorithm produces characteristic fillings which are spin, with spin structures extending the chosen spin structure on the 3-manifold.
This will be important later when we apply this construction to construct concordances.

\begin{prop} \label{prop: spin fill}
    Given any spin 3-manifold $(Y, \s)$ and null-homologous knot $K \subset Y$, there exists a spin filling $(X, \t)$ of $Y$ consisting only of 2-handles such that $K$ bounds a characteristic disk $D$ into $X$ with the property that the canonical spin characterization $\c_\t$ of $D$ in $X$ restricts to $\c_\s$ on the boundary.
\end{prop}

Before beginning the proof, we establish the following lemma.

\begin{lem} \label{lem: all knots char slice}
    For every knot $K \subset S^3$, there exists $k \in \N$ such that $K$ bounds a characteristic slice disk into $X = (\#^kS^2 \times S^2) \setminus B^4$.
\end{lem}
\begin{proof}
     Suppose first that $\Arf(K) = 0$.
     Then by Theorem \ref{thm: arf 0 band pass}, $K$ is band pass equivalent to the unknot.
     Note that we may achieve any band pass via 0-surgery on a Hopf link, which corresponds to taking a connected sum with $S^2 \times S^2$, and that in particular, each component of the Hopf link has linking number zero with $K$.
     This is illustrated in Figure \ref{fig: band pass}.

    \begin{figure}[h]
    \labellist
    \small\hair 2pt
    \pinlabel {0} at 70 10
    \pinlabel {0} at 160 100
    \pinlabel {0} at 263 120
    \pinlabel {0} at 308 140
    \pinlabel {0} at 453 130
    \pinlabel {0} at 490 140
    \endlabellist
    \centering
    \includegraphics[width=0.9\linewidth]{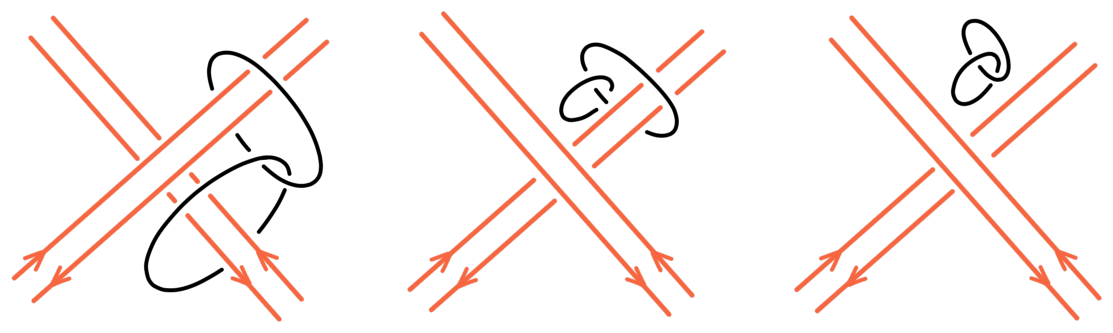}
    \caption{A band pass achieved by 0-framed 2-handles.}
    \label{fig: band pass}
    \end{figure}

    Writing $X$ as $S^3 \times I \cup \{2\text{-handles\}}$, we may obtain a disk $D$ bounded by $K \subset \partial X$ by taking the annulus $K \times I \subset X$ and capping off with the disk bounded by the unknot we are left with after all the surgeries induced by the 2-handle attachments.
    Since $K$ had linking number 0 with each 2-handle, $D \cdot x = 0$ for all $x \in H_2(X;\Z/2\Z)$.
    Finally, since $X$ is spin, $x \cdot x \equiv 0$ mod 2 for all such $x$, and so $D$ is characteristic.

    Consider now the case that $\Arf(K) =1$.
    Then again by Theorem \ref{thm: arf 0 band pass}, $K$ is band pass equivalent to the trefoil.
    We may therefore follow the construction in the previous case to obtain a properly embedded characteristic annulus $C$ embedded in a twice-punctured copy of $\#^k S^2 \times S^2$ for some $k$ such that $\partial C = K \sqcup T$ for $T$ the trefoil.
    To cap off $C$ and complete our disk, we connect sum on one more copy of $S^2 \times S^2$ to our ambient manifold, arranged so that $T$ links the 2-handles as in the diagram in Figure \ref{fig: trefoil in s2xs2}.
    
    \begin{figure}[h]
    \labellist
    \small\hair 2pt
    \pinlabel {0} at 0 60
    \pinlabel {0} at 292 60
    \endlabellist
    \centering
    \includegraphics[width=0.4\linewidth]{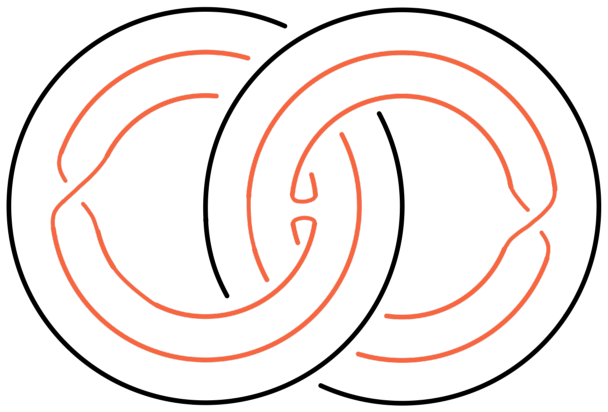}
    \caption{A trefoil in punctured $S^2 \times S^2$ which bounds a characteristic disk with relative Euler number 8.}
    \label{fig: trefoil in s2xs2}
    \end{figure}
    
    We can see that $T$ bounds a disk $D_T$ in the $(2,2)$ relative homology class of punctured $S^2 \times S^2$ with respect to the standard generators, and so $D_T \cdot x \equiv 0 \equiv x \cdot x$ mod 2 for all $x$ in $H_2((S^2 \times S^2) \setminus B^4, \Z/2\Z)$ and hence is characteristic in that copy of punctured $S^2 \times S^2$.
    Since the connected sum of characteristic surfaces is characteristic in the connected sum of their ambient manifolds, we may glue $C$ to $D_T$ along $T$ to obtain a characteristic disk bounded by $K$.
\end{proof}

Recall that the set of spin characterizations of a characteristic submanifold is a torsor over the first cohomology group of the ambient manifold. 
Thus for any knot $K$ in $S^3$, and for any characteristic disk $D$ in a simply connected manifold $X$, there is a unique spin characterization in each case.
Thus the disk constructed in the above proof has a unique spin characterization which extends the unique spin characterization on $K$.

\begin{proof}[Proof of Proposition \ref{prop: spin fill}]
    By a result of Kaplan \cite{kap79}, for any spin structure on $Y$ there always exists a spin cobordism $X_0$ between $Y$ and $S^3$ consisting of only 2-handles, and so we may interpret this cobordism as the result of spin surgery on a collection of framed curves in $Y$.

    Consider a Seifert surface $F$ for $K$ which is specified by $\c_\s$ as in subsubsection \ref{subsubsec: induced char}.
    We may assume that $F$ is disjoint from all the surgery curves via the following procedure:
    Suppose that $\gamma$ is a curve in $Y$ that we wish to do surgery along such that $\gamma \cap F$ is nonempty.
    Then for each intersection point, choose an arc in $F$ connecting the intersection point to a point on $K$. 
    Because these arcs lie in $F$ and $F$ is framed in $Y$, we may use pushoffs of these arcs as a guide to isotope $\gamma$ to remove the intersections as shown in Figure \ref{fig: 1d whitney}.

    \begin{figure}[h]
        \labellist
        \small\hair 2pt
        \pinlabel {$F$} at 110 90
        \pinlabel {$K$} at 100 20
        \pinlabel {$\gamma$} at 214 130
        \pinlabel {$\gamma$} at 10 130
        \pinlabel {$F$} at 315 90
        \pinlabel {$K$} at 310 20
        \endlabellist
        \centering
        \includegraphics[width=0.5\linewidth]{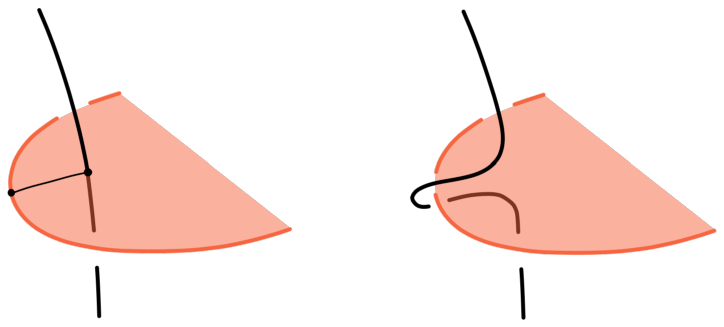}
        \caption{A 1-dimensional ``Whitney trick.''}
        \label{fig: 1d whitney}
    \end{figure}
    
    Thus we can obtain a properly embedded annulus $C$ in $X_0$ by tracing $K$ through the cobordism as in Lemma \ref{lem: all knots char slice}. 
    Specifically, writing $X$ as $Y \times I \cup \{2\text{-handles}\}$, we obtain $C$ by taking $K \times I$.
    Note that since none of the surgery curves intersect $F$, the product $F \times I$ embeds properly in $X$.
    This demonstrates that $C$ is null-homologous in $X_0$ rel. boundary, and so since $X_0$ is spin, $C$ is characteristic in $X$. 
    Furthermore, since the Seifert framing on $K$ is given by a pushoff along $F$, this framing extends over $C$.
    We conclude that $C$ is 0-framed rel. boundary.
    Finally, since $F$ is dual to $[\s - \c_\s]$, $F \times I$ is dual to $[\t - \c_\t]$ for the spin characterization $\c_\t$ induced by $\t$ of $C$ in $X$, where $\t$ is the spin structure on $X$ extending $\s$ given to us by Kaplan's algorithm. 
    Thus we may finish building the filling from $X_0$ by capping off the pair $(J, S^3)$ with a characteristic disk in some $(\#^k S^2 \times S^2) \setminus B^4$ using Lemma \ref{lem: all knots char slice}.
\end{proof}

\subsubsection{The corresponding smooth concordance problem}

Proposition \ref{prop: spin fill} shows that $\Arf(K, Y, \s)$ can be interpreted as an invariant which fully classifies null-homologous knots in spin 3-manifolds up to a very generalized notion of smooth concordance.
Although this form of concordance is so loose that it seems almost useless, we nevertheless give it a name, as it will arise naturally in our study of Montesinos annuli and periodic 2-knots.

\begin{defn} \label{def: char conc}
    We say that a triple $(K_1, Y_1, \s_1)$ is smoothly \textit{characteristically concordant} to $(K_2, Y_2, \s_2)$ if there exists a smooth spin 4-manifold $(X, \t)$ with smoothly embedded characteristic annulus $C$ cobounded by $K$ and $K'$ such that $C \cdot C \equiv 0$ mod 16 and the canonical spin characterization $\c_\t$ of $C$ in $(X, \t)$ restricts to the canonical spin characterizations $\c_\s$ and $\c_{\s'}$ of $K$ and $K'$ in $(Y, \s)$ and $(Y', \s')$, respectively.
\end{defn}

\begin{prop} \label{prop: char conc iff arf}
    The triple $(K_1, Y_1, \s_1)$ is smoothly characteristically concordant to $(K_2, Y_2, \s_2)$ if and only if $\Arf(K_1, Y_1, \s_1) = \Arf(K_2, Y_2, \s_2)$.
\end{prop}
\begin{proof}
    Suppose first that $(K_1, Y_1, \s_1)$ is smoothly characteristically concordant to $(K_2, Y_2, \s_2)$ via a spin cobordism with characteristic annulus $(C, X, \t)$. 
    For $i \in 1,2$, let $(X_i,D_i)$ be a smooth spin characteristic filling of $(K_i,Y_i,\s_i)$.
    Then the union of these disks with $C$ is a smoothly embedded characteristic sphere $S$ in the smooth spin 4-manifold given by $-X_1 \cup_{Y_1} X \cup_{Y_2} - X_2$, and hence by the Kervaire-Milnor theorem \cite{KM61} and Rokhlin's theorem \cite{roh52}, we have: 
    \[ (-D_1 \cup_{K_1} C \cup_{K_2} -D_2) \cdot (-D_1 \cup_{K_1} C \cup_{K_2} -D_2) \equiv 0 \mod 16. \]
    Thus the difference of Arf invariants of the triples $(K_i, Y_i, \s_i)$ reduces to the difference of the Rokhlin invariants $\mu(Y_i, \s_i)$, which vanishes because $(Y_1, \s_1)$ is spin cobordant to $(Y_2, \s_2)$ via $(X, \t)$.

    For the converse, suppose that $\Arf(K_1, Y_1, \s_1) = \Arf(K_2, Y_2, \s_2)$.
    By Proposition \ref{prop: spin fill}, there exist spin characteristic fillings $(X_i,D_i)$ for $(K_i, Y_i, \s_i)$.
    If we puncture these fillings at a point on $D_i$, we obtain smooth spin cobordisms between $(K_i, Y_i, \s_i)$ and $(U, S^3, \s_0)$, which may fail to be characteristic concordances only if $\Arf(K_i, Y_i, \s_i) = 1$ (in which case, the embedded annuli bounded by $K_i$ and $U$ will have self-intersection number 8 mod 16).
    In either case, we may now glue these two cobordisms together along their $(U, S^3, \s_0)$ boundaries to obtain a smooth characteristic concordance between $(K_1, Y_1, \s_1)$ and $(K_2, Y_2, \s_2)$, as the concatenation of these two annuli will have self-intersection number 0 mod 16.
\end{proof}

\subsection{Extension to arcs and additivity}

This invariant can be easily extended to properly embedded arcs in 3-manifolds with spherical boundary that are null-homologous rel. boundary.

\begin{defn} \label{def: barf arc}
    Let $Y$ be a 3-manifold with $\partial Y \cong S^2$. Let $\s$ be a spin structure on $Y$, and  let $A$ be an arc properly embedded in $Y$ such that $A$ is null-homologous in $Y$ rel. boundary.
    Consider the union of $(A,Y)$ with $(\ell, B^3)$ along their common boundary, where $\ell$ denotes the unknotted arc in $B^3$.
    Since $S^2$ and $B^3$ both have unique spin structures, we may extend $\s$ uniquely to a spin structure $\s'$ on the union.
    We define $\Arf(A, Y, \s)$ to be $\Arf(A', Y', \s')$, where $A' = A \cup \ell$ and $Y' = Y \cup B^3$.
\end{defn}

This is easily seen to be well-defined, as any two arcs that are null-homologous in their component 3-manifolds glue together to form a null-homologous knot in the union.
We also note the following additivity property, which we will need later.

\begin{lem}\label{lem: barf additivity}
    Let $K_1$ and $K_2$ be null-homologous knots in 3-manifolds $Y_1$ and $Y_2$ with spin structures $\s_1, \s_2$, respectively.
    Choose 3-balls $B_i$ embedded in $Y_i$ such that $B_i$ intersects $K_i$ in a trivial tangle, and form the connected sum of tuples $(K, Y, \s) = (K_1 \# K_2, Y_1 \#Y_2, \s_1 \# \s_2)$ along the boundary of $Y_i \setminus B_i$.
    Then:
    \[ \Arf(K,Y, \s) = \Arf(K_1, Y_1, \s_1) + \Arf(K_2, Y_2, \s_2) \mod 2.\]
\end{lem}
\begin{proof}
    For $i \in \{1,2\}$, let $(X_i,D_i)$ be a spin characteristic filling of $(K_i,Y_i,\s_i)$. 
    Choose 4-dimensional half-balls $B_i \subset X_i$ such that $B_i$ intersects $Y_i$ in a 3-ball, $K_i$ in a trivial arc contained within that 3-ball, and $D_i$ in an unknotted disk.
    If we remove $B_i$ from $X_i$, we may then glue $Y_1$ to $Y_2$ along the resulting $S^2$ boundary to obtain a closed 3-manifold $Y$, and glue the interiors of $X_1$ and $X_2$ together along the remainder of the boundary of the removed half-ball. 
    Gluing the punctured knots $K_1$ and $K_2$ and the disks $D_1$ and $D_2$ together along their boundaries resulting from the removed 4-balls, we see that $K = K_1 \# K_2$ bounds $D = D_1 \natural D_2$ into $X = X_1 \natural X_2$.
    Since $D_i$ was characteristic in $X_i$ and extended the spin characterizations on $K_i$ induced by $\s_i$, and the two disks are being glued together along a trivial arc in a ball, $D$ is characteristic in $X$ and extends the induced spin characterization of $K$ by $\s_1 \# s_2$.
    Then $D \cdot D = D_1 \cdot D_1 - D_2 \cdot D_2$ by construction, and so if we use $D$ and $X$ to compute $\Arf(K, Y, \s)$, we see that $\frac{D \cdot D}{8} = \frac{D_1 \cdot D_1}{8} + \frac{D_2 \cdot D_2}{8}$ mod 2.
\end{proof}

The above proof follows exactly the same for null-homologous arcs in 3-manifolds with spherical boundary, where the manifolds are joined via a boundary connected sum near the boundary components of each of the arcs.

This Arf invariant forms the backbone for the concordance invariants introduced in Sections \ref{sec: ann conc} and \ref{sec: equi conc}.
We will return to it then, after some constructions of Seifert manifolds.


\section{Embedded Morse theory and codimension 2 intersections} \label{sec: morse theory}

For many constructions in the following section, we will need a modified version of Morse theory to allow us to work with $n$-manifolds with boundary embedded in a cylinder of dimension $n + 2$, while keeping track of a codimension 2 submanifold.
This section is devoted to proving the foundational lemmas we will need for these constructions.

\subsection{Review of Morse theory with vertical boundary}
Morse theory with vertical boundary uses Morse functions $f$ on manifolds with boundary $M$ together with gradient-like vector fields $\xi$ which are everywhere tangent to $\partial M$.
This should not be confused with the notion of a \textit{relative} Morse function, for which $f$ is constant on connected components of $\partial M$ (e.g. on a cobordism).
These two are not mutually exclusive, and one can have a relative Morse function with vertical boundary when $M$ is a manifold with corners.
For instance, the height function $f$ on $D^2 \times I$ given by $f(x,t) = t$ is such a function.

It is important to note that we follow the convention of Kronheimer-Mrowka in \cite{km07} that the gradient-like vector field mimics the behavior of $-\nabla f$, rather than $\nabla f$.
However, we do adapt some of the language from Borodzik-N\'emethi-Ranicki in \cite{bnr16} to this convention.

Let $\varphi: M \times \R \to M$ be the flow associated to $-\nabla f$. Recall that the \textit{ascending manifold}, denoted $\calA(p)$, of a critical point $p$ is the set:
\[ \calA(p) := \{x \in M \mid \lim_{s \to \infty} \varphi(x, s) = p \} \]
Similarly, the \textit{descending manifold} of $p$, denoted $\calD(p)$ is the set:
\[ \calD(p) := \{ x \in M \mid \lim_{s\to -\infty} \varphi(x, s) = p \} \]
We introduce one more nonstandard definition which will be helpful for later sections.
\begin{defn}
    We say that an open submanifold $A$ of $M$ is \textit{fully vertical} with respect to $f$ if for all $x \in A$, $\varphi(x,s)$ is also in $A$ for all $s$ in $\R$.
\end{defn}
In other words, a fully vertical submanifold is closed under the flow of $-\nabla f$.
Note that we will use ``vertical'' as a descriptor for various manifolds throughout the paper; this is not a technical term, and is mostly used in exposition to refer to manifolds that may have critical points, but to which $\xi$ is always tangent.

\begin{defn}
    Let $M$ be a smooth $n$-manifold with boundary, and let $f$ be a Morse function on $M$.
    Let $p$ be a critical point of $f$ in $\partial M$.
    Then $p$ is said to be \textit{boundary stable} if the tangent space $T_p\calA(p)$ is contained entirely within $T_p \partial M$.
    Else, $p$ is said to be \textit{boundary unstable}.
\end{defn}

We say that the \textit{index} $\ind(p)$ of a boundary critical point is $\dim \calD(p)$.
This definition agrees with the usual index of an interior critical point.
However, note that with this definition, there can be no boundary stable critical points of index 0, nor boundary unstable critical points of index $n$.

\begin{figure}[h]
    \labellist
    \small\hair 2pt
    \pinlabel $p$ at 142 90
    \pinlabel $\calA(p)$ at 208 120
    \pinlabel $\calD(p)$ at 142 30
    \pinlabel $p$ at 430 100
    \pinlabel $\calA(p)$ at 430 160
    \pinlabel $\calD(p)$ at 500 58
    \endlabellist
    \centering
    \includegraphics[width=\columnwidth]{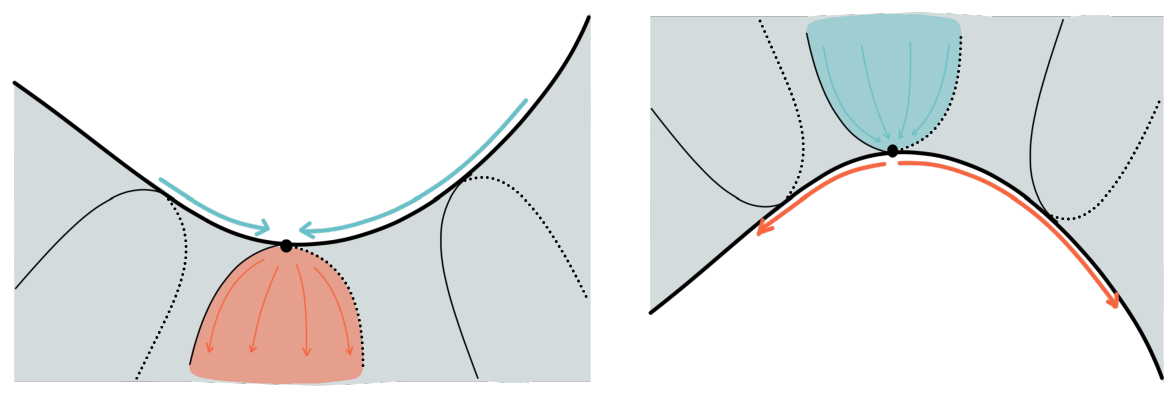}
    \caption{Left: a boundary stable critical point. Right: a boundary unstable critical point. Ascending and descending manifolds are shaded or shown as arrows according to their dimension.}
    \label{fig: boundary crit pts}
\end{figure}

\begin{defn}[\cite{bnr16}]
    An $n$-dimensional index $k$ \textit{right half-handle} $R^k$ is given by:
    \[ R ^k := D^k \times D^{n-k} \]
    When the index is clear, we simply write $R$.
    The boundary $\partial R$ is subdivided into three pieces.
    One piece is the \textit{attaching region} $A^k$, given by:
    \[ A^k \colon = \partial D^k \times D^{n-k} \]
    The remainder of the boundary is topologically a copy of $D^k \times S^{n-k-1}$.
    Decompose $S^{n-k-1}$ into two copies $D_h ^{n-k-1}$ and $D_v ^{n-k-1}$ of $D^{n-k-1}$ glued along their boundaries.
    Then the other two pieces of $\partial R^k$ we call the \textit{horizontal} and \textit{vertical} portions of the boundary, denoted $(\partial R)_h$ and $(\partial R)_v$, respectively, given by:
    \begin{align*}
        &(\partial R)_h = D^k \times D_h ^{n-k-1} \\
        &(\partial R)_v = D^k \times D_v ^{n-k-1}
    \end{align*}
\end{defn}

\begin{figure}[h] 
    \centering
    \includegraphics[width=0.75\columnwidth]{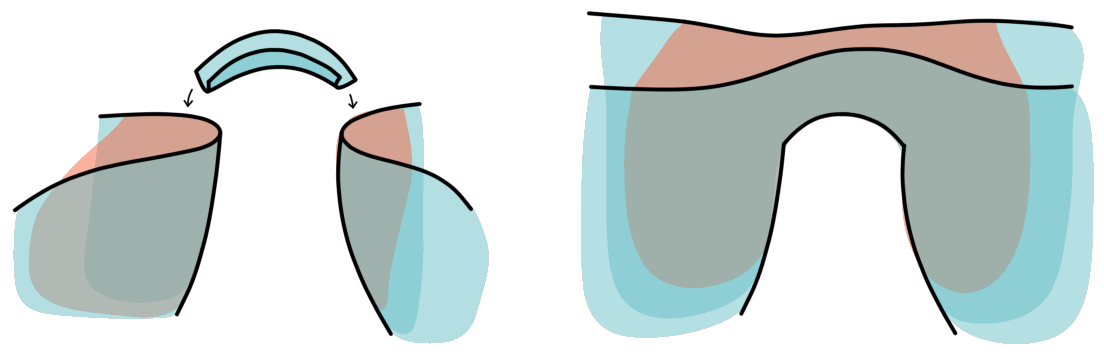}
    \caption{A 3-dimensional index 1 right half-handle attachment. Note the shading to denote the filling of the interior manifold after the boundary handle attachment.}
    \label{fig: right half-handle}
\end{figure}

This definition should be interpreted as follows. 
Suppose $M$ is an $n$-manifold with boundary and fix a relative Morse function $f$ on $M$ such that $f$ is constant on some connected component of $\partial M$.
Then an index $k$ right half-handle $R^k$ can be attached to that portion of $\partial M$ along $A^k$.
However, we do this attachment so that $(\partial R)_h$ is horizontal with respect to $f$, and $(\partial R)_v$ is vertical with respect to $f$.
Thus a right half-handle attachment changes both the topology of $\partial M$ and that of $M$ itself by attaching an index $k$ handle to each.
An example right half-handle attachment is shown in Figure \ref{fig: right half-handle}.

There is also a notion of left half-handles, but we will not use them in this paper.
See \cite{bnr16} for more details.

Borodzik-N\'emethi-Ranicki showed, in analogy with classical Morse theory, that passing a boundary unstable critical point induces a right half-handle attachment.
We will use this notion to build Seifert manifolds by hand in Section \ref{sec: seif mfds} by declaring that whenever we pass a critical point of the boundary, we attach a right half-handle to the Seifert manifold.

\subsection{Morse theory with codimension two intersections} \label{subsec: morse theory}

The goal of this subsection is to extend some of the standard results of embedded Morse theory to a setting in which we must keep track of a codimension 2 submanifold. 
Throughout, we fix the following setup:
\begin{itemize}
    \item $M$ is a smooth compact $n$-manifold, possibly with boundary or corners,
    \item $W$ is a copy of $B^{n+1} \times I$, thought of as a manifold with corners,
    \item $\psi$ is a smooth embedding of $M$ into $W$ such that there exists a decomposition of $\partial M$ into connected components and components joined by corners such that each component is either fully contained in $\partial W$, fully disjoint from $\partial W$ and has no corners, or intersects $\partial W$ only in its corners under $\psi$,
    \item $Z$ is the subset $B^{n-1} \times I$ in $W$,
    \item $A \subset M$ is a properly embedded copy of either $B^{n-2}$ or $S^{n-2}$ such that $\psi(M) \cap Z = \psi(A)$, and
    \item $f \colon W \to I$ is a relative Morse function on $W$ such that $f$ has no critical points on $W$, $f|_Z$ has no critical points, $f|_{\psi(M)}$ is Morse-Smale, and $A$ is fully vertical with respect to $f|_{\psi(M)}$ away from at most one maximum and one minimum.
\end{itemize}

\begin{figure}[h]
    \labellist
    \small\hair 2pt
    \pinlabel{$W$} at 250 100
    \pinlabel{$Z$} at 195 100
    \pinlabel{$\psi(M)$} at 120 200
    \pinlabel{$\psi(A)$} at 120 120
    \endlabellist
    \centering
    \includegraphics[width=0.5\linewidth]{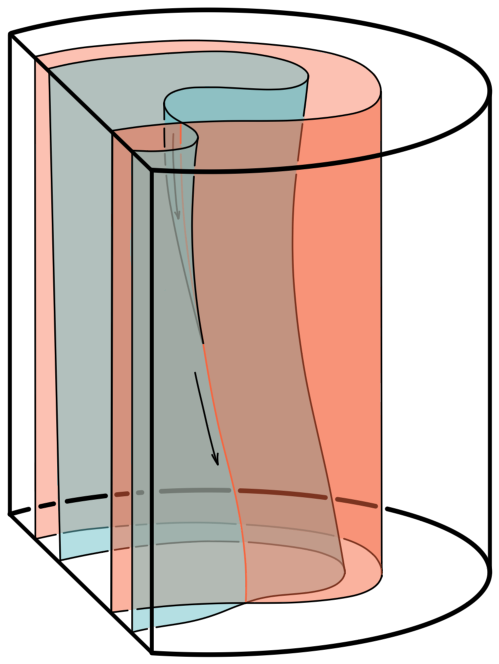}
    \caption{Schematic for the setup for Subsection \ref{subsec: morse theory} (in the case that $M$ has boundary and corners), with $f$ being projection onto the vertical axis. Note the flow lines of $-\nabla (f|_{\psi(M)})$ showing that $A$ is fully vertical with respect to $f|_{\psi(M)}$. Note that this figure is not to scale: $\psi(M)$ and $Z$ are both codimension 2 in $W$, and $A$ is codimension 2 in $M$.} 
    \label{fig: morse theory}
\end{figure}

This is a more restricted setting than appears in the literature, and so we must reprove several standard facts from standard Morse theory.
Recall the following standard definition.

\begin{defn} \label{def: omp}
    We say that $\psi(M)$ is in \textit{ordered Morse position} with respect to $f$ if:
    \begin{enumerate}
        \item the critical points of $f|_{\psi(M)}$ are isolated,
        \item the Hessian of $f|_{\psi(M)}$ at each critical point is nondegenerate, and
        \item for all critical points $p_i$ and $p_j$ of $f|_{\psi(M)}$, $f(p_i) < f(p_j)$ implies that $\ind(p_i) \leq \ind(p_j)$.
    \end{enumerate}
\end{defn}

Given our setup above, we want to be able to isotope the map $\psi$ to an embedding where $\psi(M)$ is in ordered Morse position with respect to $f$.
However, we also want this isotopy to fix a neighborhood of $A$, so that the image of $A$ is still fully vertical under the final map.
Using the convention for isotopies that $\Psi_t (x) := \Psi(x,t)$, the formal statement we want to prove is:

\begin{prop}\label{prop: omp}
    There exists an embedding $\psi' \colon M \to W$ and an isotopy $\Psi \colon M \times I \to W$ between $\psi$ and $\psi'$ such that: 
    \begin{enumerate}
        \item $\psi'(M)$ is in ordered Morse position with respect to $f$,
        \item $f|_{\psi'(M)}$ has the same number and type of critical points as $f|_{\psi(M)}$, and
        \item $\Psi_t \equiv \psi$ for all $t \in I$ on a neighborhood of $A$.
    \end{enumerate} 
\end{prop}

We prove this in two steps.
First, we show that it is possible to perturb the standard Morse function $f$ on $W$ to another Morse function which, when restricted to $M$, has ordered critical points, while avoiding perturbation on a fully vertical codimension 2 submanifold of $M$.
We then show that it is possible to achieve this same rearrangement via isotopy of $M$ inside $W$, so that $M$ is in ordered Morse position with respect to the original Morse function $f$.
Because the perturbation of the Morse function in the first step is supported outside a neighborhood of $A$, the isotopy we construct in the second step fixes $A$ pointwise, and so $A$ remains fully vertical with respect to the restriction of $f$ to the image of $M$ throughout the isotopy.

Below we have the first step.
This is essentially an adaptation of Theorem 4.1 of \cite{bp16}, and our proof follows theirs closely, but we here we are careful to keep the perturbation of the Morse function away from the fully vertical submanifold.
We also only move one critical point of the Morse function restricted to the submanifold, deviating from \cite{bp16} and other sources (e.g. \cite{mil65}) that contain this lemma that simultaneously rearrange a pair of critical points.
We use the notation for homotopies that $F_t(x) :=F(x,t)$.

\begin{lem}\label{lem: local rearrangement}
    Suppose that $f|_{\psi(M)}$ has critical points $p_k$, $p_j$ such that $k < l$ but $f(p_k) > f(p_l)$.
    Then there exists a Morse function $g$ on $W$ and a homotopy $F \colon W \times I \to \R$ such that $F_0 = f$ and $F_1 = g$, $F_t$ has no critical points on $W$ for all $t$, $g(p_k) < g(p_l)$, and with $F_t \equiv \psi$ for all $t$ on a neighborhood of $A$.
\end{lem}
\begin{proof}
    For $p \in \{p_k, p_l\}$, let $\calM(p) = \calA(p) \cup \{p\} \cup \calD(p)$.
    Then since $k < l$, a codimension counting argument shows that $\calM(p_k) \cap \calM(p_l)$ is generically empty.
    Furthermore, since $A$ was assumed to be fully vertical, both $\calM(p_k)$ and $\calM(p_l)$ are disjoint from $A$, as each of these submanifolds are fully vertical with respect to $f$ and thus closed under the flow of $-\nabla f$.

    Choose a neighborhood $N_0$ of $\calA(p_k) \cap (B^{n-1} \times \{0\})$ that is disjoint from a neighborhood of $\calA(p_l) \cap (B^{m-1} \times \{0\})$ and from a neighborhood of $A\cap (B^{m-1} \times \{0\})$.
    Let $\mu \colon B^{m-1} \times \{0\} \to \R$ be a smooth bump function supported on $N_0$ such that $\mu|_{\calA(p_k)\cap (B^{m-1} \times \{0\})} \equiv 1$.
    Then $\mu$ can be extended to a map $\nu \colon W \to \R$ by extending along the flow of a gradient-like vector field of $f$. 
    That is, if $\varphi \colon W \times \R \to W$ is the flow generated by $-\nabla f$, then:
    \[ \nu(x) = \mu \left(\lim_{s \to -\infty} \varphi(x, s)\right) \]
    In particular, $\nu$ is constant along trajectories of $f$.
    Let $N$ be the extension of $N_0$ to $W$ under this flow, so that $N$ is a neighborhood of $\calA(p_k)$ that is disjoint from neighborhoods of $\calA(p_l)$ and $A$.

    Choose $a \in [0,1]$.
    This will be our target value for our new Morse function at $p_k$.
    We construct the new Morse function by rearranging the heights of points in $W$ based on their $\nu$-value. 
    This is done by choosing a function $\Upsilon \colon [0,1] \times [0,1] \to [0,1]$ satisfying:
    \begin{enumerate}
        \item $\frac{\partial \Upsilon}{\partial r}(r,u) > 0$ for all $r$ and $u$,
        \item $\Upsilon(f(p_k), 1) = a$,
        \item $\Upsilon(r,u) = r$ for $r$ in small neighborhoods of 0 and 1 and all $u$, and also for $u$ in a small neighborhood of 0
        \item $\frac{\partial \Upsilon}{\partial r}(r,1) = 1$ for $r$ in a neighborhood of $f(p_k)$, and $\frac{\partial \Upsilon}{\partial r}(r,0) = 1$ for $r$ in a neighborhood of $f(p_l)$.
    \end{enumerate}
    Given such a function $\Upsilon$, we may define our new Morse function by:
    \[ g(w) = \Upsilon(f(w), \nu(w))\]
    Condition (1) ensures that $g$ has no critical points on $W$.
    Condition (2) is simply the statement that $g(p_k) = a$.
    Condition (3) requires that $f$ and $g$ agree on the horizontal boundary of $\partial W$ and away from $N$.
    Finally, condition (4) ensures that $g$ is essentially a linear shift of $f$ near the critical points $p_k$ and $p_l$, which preserves their type (interior vs. boundary) and index.

    We see that since $g$ and $f$ agree away from $N$, they in particular agree on a neighborhood of $A$.
    Through a continuous family of choices for values of $a$ terminating with $a < f(p_l)$, we find a family of Morse functions with no critical points on $W$ interpolating between $f$ and a Morse function $g$ with $g(p_k) < g(p_l)$, and with $A$ remaining fully vertical throughout.
\end{proof}

We thank Mark Powell for his help with the following lemma, which details the second step of the outlined procedure.

\begin{lem} \label{lem: morse fn to isotopy}
    Suppose that $f$ and $g$ are two relative Morse functions on $W$ that each have no critical points on $W$ and are homotopic through submersions of $W$.
    Then there exists an isotopy $\Psi \colon M \times [0,1] \to W$ and an orientation-preserving diffeomorphism $h \colon \R \to \R$ such that $\Psi(x,0) = \psi(x)$ and $\Psi(x,1)$ satisfies:
    \[ f\circ \Psi(x,1) = h\circ g\circ \psi(x) \]
\end{lem}
Note that, as in Lemma \ref{lem: local rearrangement}, $f$ and $g$ may have critical points when restricted to $\psi(M)$, but cannot have any on $W$ itself.
We remind the reader of the definition of a stable map to aid in the proof.
For more on stable maps, the author recommends \cite{gg73}.

\begin{defn}
    Let $f \in C^\infty(X,Y)$ for $X$ and $Y$ smooth manifolds. 
    Then $f' \in C^\infty(X,Y)$ is said to be \textit{equivalent} (in the sense of stable maps) to $f$ if there exist orientation-preserving diffeomorphisms $\phi \colon X \to X$ and $\eta \colon Y \to Y$ such that the following diagram commutes:
    \[\xymatrix{
    X \ar[r]^-{f} \ar[d]_{\phi} & Y \ar[d]^{\eta} \\
    X \ar[r]_-{f'} & Y\\} \]
    A map $f$ is said to be \textit{stable} if there exists a neighborhood $N(f) \subset C^\infty(X,Y)$ such that all $f' \in N(f)$ are equivalent to $f$.
    \end{defn}

\begin{proof}[Proof of Lemma \ref{lem: morse fn to isotopy}]
    Let $F \colon W \times [0,1] \to \R$ be a homotopy through relative Morse functions with no critical points such that $F(w,0) = f(w)$ and $F(w,1) = g(w)$.
    Then $F(w,t)$ is Morse for all $t$.
    Morse functions are stable so long as they do not take the same value on any two critical points \cite[Ch. III Prop. 2.2]{gg73}, so for each $t$, there exists a neighborhood $N(F(-,t))$ in which all functions are equivalent to $F(-,t)$.
    We therefore have that all the Morse functions $F(-,t)$ are equivalent to each other, and so for each $t$ we have maps $\phi_t$ and $\eta_t$ such that the following diagram commutes:
    \[\xymatrix{
    W \ar[r]^-{f} \ar[d]_{\phi_t} & \R \ar[d]^{\eta_t} \\
    W \ar[r]_-{F(-,t)} & \R\\} \]
    So if we let $\Psi(x,t) = \phi_t ^{-1}(\psi(x))$, then our diagram becomes:
    \[\xymatrix{
    & W \ar[r]^-{f} \ar[d]_{\phi_t} & \R \ar[d]^{\eta_t} \\
    M \ar[r]_-{\psi} \ar[ur]^-{\Psi(-,t)} & W \ar[r]_-{F(-,t)} & \R\\} 
    \]
    And so we see that:
    \[ f\circ \Psi(x,1) = \eta_1^{-1}\circ F(\psi(x),1) = \eta_1 ^{-1}\circ g \circ \psi(x) \]
    This completes the proof.
\end{proof}

For our purposes, we need not worry about the function $\eta_1 ^{-1}$ when finding an isotopy interpolating between $f$ and $g$, as orientation preserving diffeomorphisms of $\R$ are all isotopic to the identity and merely contribute a rescaling of our Morse function that we can undo later.

\begin{rem}\label{rem: support}
    From the proof of stability for Morse functions with distinct critical values (see \cite[Ch. III Prop. 2.2]{gg73}), one sees that the functions $\phi_t$ are constructed so that they are supported on neighborhoods of the subsets of $W$ on which $F(-,t)$ differs from $f$.
    Thus the image of the support of the isotopy $\Psi$ can be taken to be contained in a neighborhood of the support of $f - g$.
\end{rem}

We can now prove the existence of isotopies to ordered Morse position in the complement of a codimension two intersection.

\begin{proof}[Proof of Proposition \ref{prop: omp}] 
    Recall that we begin with a Morse function $f$ on $W$, and that $\psi(A) = \psi(M) \cap Z$.
    By assumption, $f|_{\psi(A)}$ has at most two critical points; denote these $p_+$ and $p_-$ according to whether they are maxima or minima.
    If $p_+$ exists, choose a neighborhood of $p_+$ and perturb $\psi$ to isotope it vertically along a path tangent to $\nabla f$ in $Z$ such that $f(p_+) > f(x)$ for all $x \in \psi(M)$, using a bump function to smooth the corners.
    Similarly, if $p_-$ exists, perturb $\psi$ to isotope a neighborhood of $p_-$ vertically down along a path tangent to $-\nabla f$ so that $f(p_-) < f(x)$ for all $x$ in $\psi(M)$ and smooth the corners with a bump function.
    This procedure guarantees that the lowermost index 0 critical point and the uppermost index $n$ critical point of $f|_{\psi(M)}$ are both in $\psi(A)$.

    Let $a_\pm = f(p_\pm)$, and let $a'_\pm = a_\pm \mp \epsilon$ for $\epsilon$ small enough that the only two critical points of $f|_{\psi(M)}$ outside of $V:= [a_-', a_+']$ are $a_-$ and $a_+$.
    Then $M \cap V$ is a cobordism between two unknotted spheres (or possibly disks, if $\partial M \neq \emptyset$) embedded in $B^{n+1} \times [a_-',a_+']$.
    
    Since $f$ has no critical points on $W$, we may use Lemma \ref{lem: local rearrangement} to successively perturb $f$ to obtain Morse functions $f_i$ that rearrange pairs of critical points in the complement of $A$.
    The perturbations are supported on neighborhoods of descending manifolds of lower index critical points, and each $f_i$ is related to $f_{i+1}$ by a homotopy through submersions.
    Hence by Lemma \ref{lem: morse fn to isotopy}, we may find an isotopy $\Psi \colon M \times [0,1] \to W$ realizing this rearrangement, with $\Psi_0 = \psi$. 
    By the Isotopy Extension Theorem, this can then be extended to an ambient isotopy which, by abuse of notation, we also call $\Psi$ .
    As we noted in Remark \ref{rem: support}, since each perturbation of the Morse function was disjoint from $Z$, $\Psi$ can be taken such that $\Psi_t \equiv \psi$ for all $t$ on a neighborhood of $A$. 
    This completes the proof.
\end{proof}

\subsection{The equivariant case}

It is worth saying a few words about equivariant Morse theory for finite cyclic group actions, and detailing why we may content ourselves with the lemmas proved in the previous subsection when working with periodic 2-knots.
It should be noted that in general, equivariant differential topology is much more restrictive than the non-equivariant setting.
Transversality results often fail equivariantly.
However, with sufficient conditions on the group action, the fixed set, and the manifold acted on, many standard results in Morse theory can be recovered in the equivariant setting.
For instance, in the context of equivariant 4-manifold topology, Meier and Scott were able to prove an equivariant version of the Laudenbach-Po\'enaru theorem \cite{ms25a}, which enabled them to study equivariant trisections and bridge trisected surfaces \cite{ms25b}.
We will not use these results explicitly, but we mention them due to their close relationship to our study of equivariant 2-knots.

Firstly, note that for our purposes and in much of the literature of equivariant Morse theory, the action of $\Z/d\Z$ on the codomain $\R$ is assumed to be trivial.
Thus when we say that a Morse function $f$ is \textit{equivariant}, we mean that $f \circ \rho = f$ where $\rho$ is the action on the domain.
This was the setting studied by by Wasserman in \cite{was69}, building on work of Bott \cite{bot54} and Palais \cite{pal63}.
Also worth mentioning is work of Bao and Lawson \cite{bl24} in which they gave an extension of Wasserman's result that shows genericity of Morse-Smale functions among pairs of equivariant functions and equivariant metrics.

The key fact we take from Wasserman's work is that for compact Lie groups $G$, the correspondence between critical points of a Morse function and handle attachment carries over to the $G$-equivariant setting by defining a critical $G$-orbit.
In the case of our finite order rotation actions $\rho$ with $\Fix(\rho)$ either a sphere or a ball, equivariant handles always appear as either a single handle that is fixed setwise by $\rho$, or as a $d$-tuple of handles that are permuted by $\rho$.
Therefore, we may work with handle decompositions of simply periodic manifolds simply by working with their quotients and then lifting to the cover, so long as those handle decompositions are well-behaved relative to the image of the fixed set.
Furthermore, it is immediate that $\Fix(\rho)$ is fully vertical with respect to $f$ away from its maximum and minimum.
This, combined with Schwarz's theorem on lifting isotopies, allows us to port all the results from the previous section into the $\Z/d\Z$-equivariant setting.
For instance, we may define and prove the existence of an ordered equivariant Morse position for simply periodic manifolds which is completely analogous to Definition \ref{def: omp}.
Again, we fix the following notation for this subsection:
\begin{itemize}
    \item $\tilde M$ is a smooth, compact, simply-$d$-periodic $n$-manifold, possibly with boundary or corners,
    \item $\tilde W$ is a copy of $\tilde B^{n+1} \times I$, thought of as a manifold with corners, equipped with the linear order $d$ rotation $\rho$ fixing $\hat B^{n-1} \times I$,
    \item $\tilde \psi$ is a smooth, equivariant embedding of $\tilde M$ into $\tilde W$ such that there exists a decomposition of $\partial \tilde M$ into connected components and components joined by corners such that each component is either fully contained in $\partial \tilde W$, fully disjoint from $\partial \tilde W$ and has no corners, or intersects $\partial \tilde W$ only in its corners under $\tilde  \psi$,
    \item $\hat A \subset \tilde M$ is a properly embedded copy of either $B^{n-2}$ or $S^{n-2}$ such that $\tilde \psi(\tilde M) \cap \Fix(\rho) = \tilde \psi(\hat A)$ (and therefore $\hat A$ is the fixed set of $\tilde M$ under its periodic action), and
    \item $\tilde f \colon \tilde W \to I$ is a relative Morse function on $\tilde W$ such that $\tilde f$ has no critical points on $\tilde W$, $\tilde f|_{\Fix(\rho)}$ has no critical points, $\tilde f|_{ \tilde \psi(\tilde M)}$ is Morse-Smale, and $\hat A$ is fully vertical with respect to $\tilde f|_{\tilde \psi(\tilde M)}$ away from at most one maximum and one minimum.
\end{itemize}

\begin{defn} 
    We say that $\tilde M$ is in \textit{ordered equivariant Morse position} with respect to $\tilde f$ if:
    \begin{enumerate}
        \item the critical points of $\tilde f|_{\tilde M}$ are isolated,
        \item the Hessian of $\tilde f|_{\tilde M}$ at each critical point is nondegenerate,
        \item every critical point of $\tilde f|_{\tilde M}$ is either fixed by $\rho$ or sent to another critical point of the same index, and
        \item for all critical points $p_i$ and $p_j$ of $\tilde f$, $\tilde f(p_i) < \tilde f(p_j)$ implies that $\ind(p_i) \leq \ind(p_j)$.
    \end{enumerate}
\end{defn}

The following is the equivariant counterpart to Proposition \ref{prop: omp}.

\begin{prop}\label{prop: oemp}
    There exists an embedding $\tilde \psi' \colon \tilde M \to \tilde W$ and an equivariant isotopy $\tilde \Psi \colon \tilde M \times I \to \tilde W$ between $\tilde \psi$ and $\tilde \psi'$ such that: 
    \begin{enumerate}
        \item $\tilde \psi'(\tilde M)$ is in ordered Morse position with respect to $\tilde f$,
        \item $\tilde f|_{\tilde \psi'(\tilde M)}$ has the same number and type of critical points as $\tilde f|_{\tilde \psi(\tilde M)}$, and
        \item $\tilde \Psi_t \equiv \tilde \psi$ for all $t \in I$ on a neighborhood of $\hat A$.
    \end{enumerate} 
\end{prop}

\begin{proof}
    Since $\tilde f$ and $\tilde \psi$ are both equivariant, we may take the quotient our whole setup by $\rho$.
    We then obtain an embedding $\psi \colon M \to W$ such $\psi(M) \cap \Branch(\rho)$ is just the image $\psi(A)$ of $\hat A$ after taking the quotient of $\tilde M$ by its periodic action.
    We are then in precisely the setup for Proposition \ref{prop: omp}, and so there exists an isotopy $\Psi$ of $M$ in $W$ rearranging $M$ to ordered Morse position with respect to the quotient Morse function $f$.
    Finally, we can then take the cyclic branched covers of $M$ and $W$ over their respective fixed sets and take $\tilde \Psi$ to be the induced isotopy upstairs, which lifts by Schwarz's resolution of the Isotopy Lifting Conjecture \cite{sch80} because $\Psi$ is supported away from $A$.
\end{proof}


\section{Seifert 3- and 4-manifolds} \label{sec: seif mfds}

For the discussion of concordance annuli and later of periodic 2-knots in Section \ref{sec: equi conc}, we will need the existence of certain smooth Seifert manifolds, both in ambient dimensions 4 and 5. 
The Seifert 3-manifolds we construct are for spheres which represent annuli, that is, we build a Seifert 3-manifold while keeping track of its intersection with the removed sphere.

\subsection{Conventions for the constructions}

To aid in the following constructions, we make some parametrization conventions.
Let $U$ be a 5-dimensional half-ball in $B^5$ intersecting the boundary in a 4-dimensional ball.
We parametrize $U$ by the coordinates:
\[ U = \{ (x, y, z, w, t)\mid x^2 + y^2 + z^2 + w^2 + t^2 \leq C, t \geq 0 \} \]
for some $C \gg1$.
When we are doing constructions explicitly in $U \cap S^4$, we will refer to coordinates on $U$ simply as $(x,y,z,w)$.
We also implicitly choose the standard metric on $U$, scaled appropriately to our choice of coordinates, which we use in our discussion of Morse theory.
Submanifolds inside $U$ inherit this metric, and we will frequently differentiate between the gradient of a Morse function on the ambient space and the gradient of its restriction to the submanifold.

Finally, we use the following convention: if $f$ is a Morse function on a manifold $M$, then:
\[M_{c} := f^{-1}(c)\]
for $c \in \R$.
We also declare:
\[M_{*c} := f^{-1}(\{a \in \R \mid a * c\})\]
where $*$ is an element of $\{<, \leq, >, \geq \}$.

\subsection{Seifert 3-manifolds intersecting a plane}

\begin{prop} \label{prop: seif 3-mfds}
    Let $S$ be a smooth 2-knot in $S^4$, fix $R$ a smooth unknotted 2-sphere in $S^4$, and suppose that $S \cap R$ consists of exactly two points $p_+$ and $p_-$.
    Then $S$ bounds a smoothly embedded Seifert 3-manifold $Y$ in $S^4$ such that $Y \cap R$ is an arc $A$ bounded by $p_+$ and $p_-$ which is null-homologous in $Y$ rel. boundary.
\end{prop}
\begin{proof}
    Fix a Morse function $f$ on $S^4$ such that its restriction to a ball around $S$ has exactly two critical points.
    By Proposition \ref{prop: omp}, we may assume that $S$ is in ordered Morse position with respect to $f$.
    As noted in the previous subsection, we use coordinates $(x,y,z,w)$ such that $f$ is projection onto the $w$-axis and $R$ appears as the $(z,w)$-subspace.

    Recall that by $S_w$, we mean the preimage of $w$ under $f$ intersected with $S$, thought of as living in the whole 3-dimensional preimage of $w$.
    We call this preimage a $w$-slice.
    We may construct $Y$ by specifying $Y_w$, defined analogously, for each $S_w$.
    In particular, for regular values $w$ of $f$, $S_w$ is a link which does not intersect the $z$-axis, and $Y_w$ is a Seifert surface for that link intersecting the $z$-axis in a single point.
    A smoothly varying family of such surfaces then stack together to form $Y$ with the desired fixed set.
    We will regularly switch between this perspective and explicit handle attachments to $Y$ throughout the construction.

    \bigskip

    \noindent $\pmb{f = 0.}$
    Since $S$ is in ordered Morse position, we may consider all of its 0-handles first. 
    These will be a collection of unlinked, unknotted disks.
    For each of these disks, we begin our construction of $Y$ by attaching a 3-dimensional index 0 right half-handle.
    Since all these attachments occurred in $f^{-1}(0,\epsilon)$, we see that $S_{\nicefrac{1}{2}}$ is a collection of unlinked, unknotted circles such that one of them links the $z$-axis, and $Y_{\nicefrac{1}{2}}$ appears as a collection of standardly embedded disks bounded by each component of $S_{\nicefrac{1}{2}}$.

\bigskip

    \noindent $\pmb{f = \frac{1}{2}.}$
    Now for each index 1 critical point, we consider the corresponding descending manifold.
    These consist of collections of arcs passing through the index 1 critical points, possibly linking components of $S_w$ for $w$ below the critical value, and eventually intersecting those unknotted components.
    Up to vertical isotopy, these intersections can be chosen to occur in $f^{-1}(\frac{1}{2}, \frac{1}{2}+\epsilon)$.

    Because the descending manifolds of the index 1 critical points of $S$ may link the components of $S_w$ for $\frac{1}{2}+ \epsilon < w < 1$, we cannot simply continue filling in the 0-handles $S$ with disks past $w=\frac{1}{2}$, as they would then intersect the descending manifolds. 
    To avoid this, we first perturb $S$ to introduce new intersections of the descending manifolds with disks bounded by $S_{\nicefrac{1}{2}}$ as necessary so that each descending arc intersects the disks in pairs of oppositely signed points.
    We then attach 3-dimensional \textit{interior} 1-handles to $Y$ whose corresponding interior critical points have critical values in $(\frac{1}{2}-\epsilon, \frac{1}{2})$ and whose attaching spheres are vertical translates of the intersections of the index 1 boundary trajectories with the disks bounded by $S_{\nicefrac{1}{2}}$. 
    In each $w$-slice for $w$ above $\frac{1}{2}$, we therefore see tubes attached to the disks bounded by $S_w$ so that the descending manifolds of the index 1 boundary critical points are engulfed by these tubes and thus do not intersect our Seifert surface in those cross-sections.

    \begin{figure}[h]
    \centering
    \includegraphics[width=0.9\columnwidth]{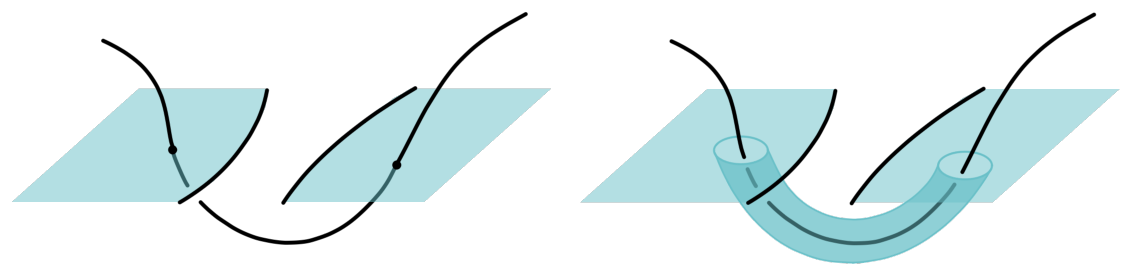}
        \caption{Tubes attached to $Y_{\nicefrac{1}{2}}$ removing intersections with descending manifolds as a result of interior 1-handle attachments.}
        \label{fig: tubing to avoid bands}
    \end{figure}

    \bigskip

    \noindent $\pmb{f = 1.}$
    We may now consider each index 1 critical point of $S$ as a boundary unstable index critical 1 point of $Y$, and thus for each 1-handle of $S$, we attach one right half-handle to $Y$.
    After performing the above procedure for all 1-handles in $S$, the cross-section $S_{\nicefrac{3}{2}}$ will consist again of a collection of unknotted components, one of which is a meridian for the $z$-axis.
    The corresponding surface $Y_{\nicefrac{3}{2}}$ will be a possibly disconnected surface of unknown genus with boundary $S_{\nicefrac{3}{2}}$.

    \bigskip

    \noindent $\pmb{f = \frac{3}{2}.}$
    At this point, $S_{>\nicefrac{3}{2}}$ consists only of disks, but $Y_{\nicefrac{3}{2}}$ is a potentially complicated surface.
    We would like to cap it off with 3-handles, but to do this, we must first surger it to a collection of disks.
    Therefore we now proceed to find a collection of compressing curves on $Y_{\nicefrac{3}{2}}$ such that the compressions do not introduce any new intersections with $R$.

    Let $\{J_i\}_{i \in I}$ be the set of connected components of $S_{\nicefrac{3}{2}}$, with $J_0$ the component that links the $z$-axis.
    Then the components $\{J_i\}_{i \in I}$ bound a pairwise disjoint collection of maximum disks $\{D_i\}_{i \in I}$ of $S$ such that only $D_0$ intersects $R$, and it does so transversely and exactly once.
    We may therefore isotope copies of these disks into $f^{-1}(\frac{3}{2})$ by flowing along $-\nabla f$.
    Call these isotoped disks $ D_i '$.
    Then $\epsilon$-neighborhoods of the disks $D_i '$ for $\epsilon$ sufficiently small have boundaries that form splitting spheres $V_i$ separating each component of $S_{\nicefrac{3}{2}}$ from each other, and separating all $J_i$ except $J_0$ from the $z$-axis in the $\frac{3}{2}$-slice.
    
    We can now make a standard ``innermost circles'' argument to find compressing disks for $Y_{\nicefrac{3}{2}}$.
    Consider the intersections of each splitting sphere $V_i$ with $Y_{\nicefrac{3}{2}}$.
    These intersections form nested concentric circles on each sphere $V_i$.
    By compressing along these concentric disks starting with the innermost circle of intersection and working outward, we can systematically compress $Y_{\nicefrac{3}{2}}$ without introducing new intersections with the splitting spheres, as shown in Figure \ref{fig: compressing Y}.
    Since these spheres are disjoint from the $z$-axis, and lie entirely in a single $w$-slice, the compressing disks are also disjoint from $R$, and so compressing along them does not introduce any new intersections with $R$. 
    Actually performing the compressions is done by attaching interior 2-handles to $Y_{\leq \nicefrac{3}{2}}$, where the core disks of the handles are taken to be the compressing disks themselves.
    We declare that these handle attachments occur within $f^{-1}(\frac{3}{2}, \frac{3}{2} + \epsilon)$. 
    
    \begin{figure}[h]
        \labellist
        \small\hair 2pt
        \pinlabel $J_i$ at 55 75
        \pinlabel {$D_i '$} at 55 105
        \pinlabel{$V_i$} at 200 50
        \endlabellist
        \centering
        \includegraphics[width=\linewidth]{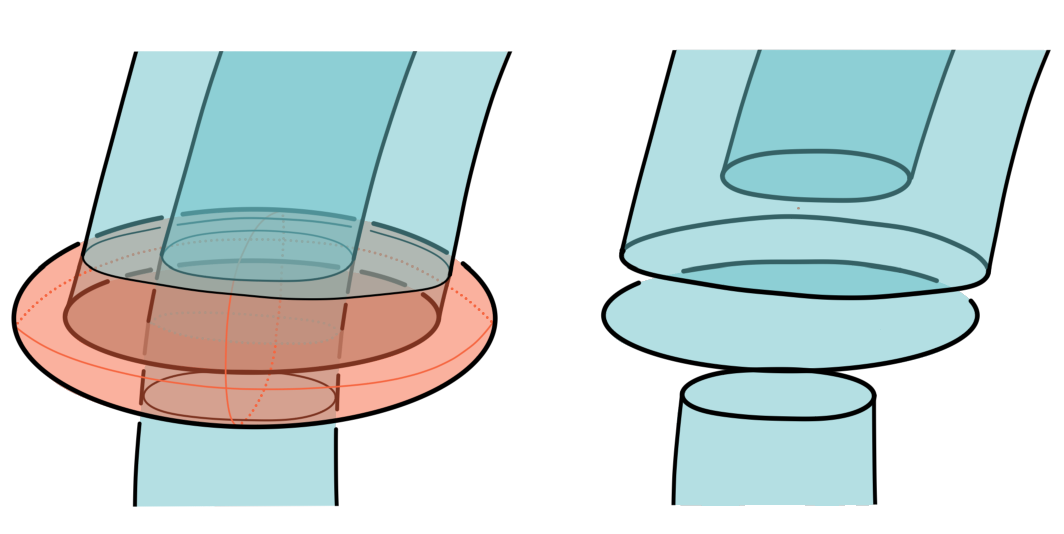}
        \caption{``Innermost circles'' compression of $Y_{\nicefrac{3}{2}}$ using compressing disks obtained via the boundary sphere $V_i$ of a tubular neighborhood of one of the disks $D'_i$.}
        \label{fig: compressing Y}
    \end{figure}
    
    \bigskip
    \noindent $\pmb{f = 2.}$
    After all such compressions are performed, $S_{\nicefrac{3}{2} + \epsilon}$ is unchanged, but $Y_{\nicefrac{3}{2} + \epsilon}$ now consists of a collection of disjoint disks and closed components, still intersecting the $z$-axis in a single point contained inside a single disk.
    We further compress the closed components until they are spheres, and then fill them in with interior 3-handles.
    This also does not change $S_{\nicefrac{3}{2}+\epsilon}$, but surgers $Y_{\nicefrac{3}{2} + \epsilon}$ to a collection of disks that are isotopic to the maximum disks of $S$.
    Thus we may fill in the remainder of $S$ by declaring $Y_{\geq \nicefrac{3}{2} + \epsilon}$ to be trace of this isotopy.
    This gives $Y$, and its intersection with $R$ is precisely its intersection with the $w$-axis, which is a single arc as desired.
    
    \bigskip
    
    It remains to show that the intersection $A$ of $Y$ with $R$ is null-homologous in $Y$ rel. boundary.
    Indeed, $A$ does not intersect the belt spheres of any 1-handles of $Y$, as all of the 1-handles were chosen to miss $R$.
    Thus $A$ represents the 0 class in the cellular homology rel. boundary induced by the handle decomposition.
\end{proof}

\subsection{Seifert 4-manifolds intersecting a 3-ball}

We now prove a similar statement to Proposition \ref{prop: seif 3-mfds}, but one dimension up.
Note that for this construction, there are extra Morse-theoretic conditions that were not present in the previous.
It is unknown at this time whether these can be removed, and this question is closely tied to the discussion in Remark \ref{rem: equi conc def}.

\begin{prop} \label{prop: seif 4-mfds}
    Fix $Q$ a smoothly unknotted 3-sphere in $S^5$.
    Let $Y$ be a closed 3-manifold smoothly embedded in $S^5$, and suppose that $Y \cap Q$ is a null-homologous knot $K$ in $Y$.
    Suppose that there exists a Morse function $f$ on $S^5$ such that $Y$ is in ordered Morse position with respect to $f$, $f|_K$ has exactly two critical points, and away from these two critical points, $K$ is fully vertical with respect to $f|_Y$.
    Then $Y$ bounds a smooth Seifert 4-manifold $X$ smoothly embedded in $S^5$ such that $X \cap Q$ is a smoothly embedded disk $D$ bounded by $K$ which is null-homologous in $X$ rel. boundary.
\end{prop}
\begin{proof}
    We follow roughly the same strategy as in the proof of Proposition \ref{prop: seif 3-mfds}. 
    That is, we build $X$ handle by handle according to a handle decomposition of $Y$ coming from $f$.
    We do this, as before, by finding fillings for slices of $Y$, which themselves are surfaces.
    To better visualize these surfaces, we parametrize a neighborhood of $Y$ so that $Q$ appears locally as the $(z,w,t)$-subspace, and so that $f$ is just projection to the $t$ coordinate.
    Note that the full verticality condition on $K$ ensures that the lowest minimum (maximum) of $Y$ is necessarily also the minimum (resp. maximum) of $K$, and so every $t$-slice has at least one intersection point with $Q$ in it.
    The assumption that $K$ is fully vertical with respect to $f|_Y$ away from two critical points implies that no $t$-slice of $Y$ intersects $Q$ in more than two points.
    Thus in each $t$-slice apart from those containing the very lowest minimum and the very highest maximum, $Y$ appears as a surface in 4-space intersecting the $(z,w)$-subspace in two points.
    Our goal is therefore to construct smoothly varying 3-dimensional fillings $X_t$ of these surfaces $Y_t$ such that each of these slices intersects the $(z,w)$-subspace in a single smoothly varying arc.
    These arcs stacked on top of each other will then form a disk $D$, which will be $X \cap Q$ as desired.
    
    \bigskip
    
    \noindent $\pmb{f = 0.}$
    Since $Y$ is in ordered Morse position, we can assume that all of its 0-handles occur in $0 \leq t < \epsilon$, with the lowest minimum intersecting $Q$ in an unknotted arc.
    We may begin constructing $X$ by filling in these 0-handles with 4-dimensional index 0 right half-handles, such that the aforementioned arc bounds a disk given by a straight line in the $(z,w)$-plane connecting each intersection point at each level.
    This completes $X$ and $D$ up to $t= \frac{1}{2}$.
    
    \bigskip
    
    \noindent $\pmb{f = \frac{1}{2}.}$
    At this point, $X_{\nicefrac{1}{2}}$ is a collection of 3-balls in ambient 4-space.
    As in the previous case, we want to fill the 1-handles of $Y$ with 4-dimensional right half-handles. 
    However, the trajectories from the index 1 critical points of $Y$ may link slices of the 0-handles of $Y$.
    As before, we consider these trajectories intersecting the 0-handles of $Y$ in $f^{-1}(\frac{1}{2}, \frac{1}{2} + \epsilon)$.
    To avoid these intersections, we preemptively attach interior 1-handles to $X$ at $t= \frac{1}{2} - \epsilon$ along the trajectories to create 3-dimensional tubes for these trajectories to pass through.
    Since the trajectories from the 1-handles of $Y$ do not intersect $Q$, and all interior 1-handles are attached parallel to these trajectories, we can assume that none of the 1-handles of $X$ intersect $Q$.
    Thus we may continue constructing $D$ by taking the straight line in the $(z,w)$-plane connecting $Y_t \cap Q$ in each $t$-slice.
    
    \bigskip
    
    \noindent $\pmb{f = 1.}$
    After these interior 1-handles are attached, we may proceed to attach 4-dimensional index 1 right half-handles to fill each 1-handle of $Y$.
    Again, since none of the 1-handles of $Y$ intersect $Q$, these half-handles do not either, and we can continue constructing $D$ by taking the straight line in the $(z,w)$-plane connecting the two points of $Y_t \cap Q$ in each $t$-slice.

    
    \bigskip
    
    \noindent $\pmb{f = \frac{3}{2}.}$
    Following the same general strategy, we now examine the descending manifolds of the 3-dimensional index 2 critical points at $t = 2$. 
    Because $Y \cap Q$ is precisely $K$, we know none of these descending manifolds intersect $Q$, as all such intersections are already accounted for by our assumption that $K$ has precisely two critical points with respect to the ambient Morse function.
    Thus in theory we may safely attach right half-handles to $X$ for each 2-handle of $Y$ without fear of introducing new components of $X \cap Q$.
    However, the index 2 descending manifolds may intersect interior 1-handles of $X$ that we have already attached. 
    As before, suppose that all of these intersections occur within the preimage $f^{-1}(\frac{3}{2}, \frac{3}{2} + \epsilon)$.
    
    To avoid these intersections, consider an index 2 critical point of $Y$. 
    Note that the generic intersection between a 2-handle core and a 4-dimensional 1-handle in $B^5$ is a disk in the 2-handle core.
    Thus we can use this disk as a guide to preemptively attach 4-dimensional interior 2-handles to compress the 4-dimensional 1-handles to avoid intersecting the 3-dimensional 2-handle.
    This is completely analogous to the preemptive 1-handles we attached to make way for the 1-handles of $Y$.
    Note further that because the disks we use as guides for attaching our 4-dimensional 2-handles are already embedded in $B^5$ as the cores of existing 3-dimensional 2-handles, they have trivial normal bundles.
    Thus from the perspective of $X$, the compressions are achieved by attaching interior 2-handles, with framings given by a choice of trivial disk sub-bundle of the normal bundle.
    Furthermore, since none of the 2-handles of $Y$ intersect $Q$, we may arrange so that the parallel 2-handles of $X$ do not intersect $Q$ either.
    
    \bigskip
    
    \noindent $\pmb{f = 2.}$
    Once we have performed these compressions at $t=\frac{3}{2}-\epsilon$, we may safely proceed upwards in $t$ and fill the 3-dimensional 2-handles of $Y$ with 4-dimensional index 2 half-handles at $t=2$.

    
    \bigskip
    
    \noindent $\pmb{f = \frac{5}{2}.}$
    At this point, for $t \in (2+ \epsilon, \frac{5}{2} - \epsilon)$, $Y_t$ appears as an unlink of unknotted spheres intersecting the $(z,w)$-plane in exactly two points on a single component, and $X_t$ appears as a connected 3-manifold with boundary $Y_t$ intersecting the $(z,w)$-plane in a single straight line.
    Now every component of $Y_{\nicefrac{5}{2}-\epsilon}$ bounds a 3-ball embedded in $B^5_{\geq 2+\epsilon}$ given by the descending manifolds of each maximum of $f|_Y$.
    As in the 3-manifold case, we may isotope these 3-ball caps straight down into the $t=\frac{5}{2}-\epsilon$-slice, using their boundary unknotted 2-spheres as guides near the boundary.
    
    As in the lower-dimensional case, the boundaries of tubular neighborhoods of each of these 3-balls then form splitting spheres for $Y_{\nicefrac{5}{2}-\epsilon}$.
    These splitting spheres intersect $X_{\nicefrac{5}{2}-\epsilon}$ in a collection of 3-balls, which we may then use as compressing 3-balls for $X_{\nicefrac{5}{2}-\epsilon}$.
    We can then compress $X$ along these 3-balls by attaching interior 3-handles to $X_{\nicefrac{5}{2}}$ and filling in closed components.
    
    \bigskip
    
    \noindent $\pmb{f = 3.}$
    In fact, because the splitting spheres we found were the boundaries of small neighborhoods of balls bounded by components of $Y_{\nicefrac{5}{2}}$, each component of $X_{\nicefrac{5}{2}}$ with boundary is itself a 3-ball.
    We may subsequently fill the remaining closed components with equivariant 4-manifold fillings.
    Finally, we trace the isotopy of the 3-ball caps up to the maxima of $Y$ to find 4-balls that fill these maxima, with the central 4-ball intersecting $Q$ in a disk.
    This completes the construction of $X$.
    
    We now consider the disk $D$ that we have constructed as $\Fix(\rho) \cap X$.
    As in the case of $A$ in the proof of Proposition \ref{prop: equi seif 3-mfds}, $D$ does not intersect the belt spheres of any handles used in the construction of $X$, and thus represents the 0 class in the relative homology of $X$.
\end{proof}

\subsection{Computing the Arf invariant using Seifert 4-manifolds}

The purpose of Proposition \ref{prop: seif 4-mfds} is to show that the Arf invariant defined in Section \ref{sec: barf} can be computed using a Seifert 4-manifold:

\begin{prop} \label{prop: seif 4-mfds give barf 0}
    Let $Y, Q$, and $K$ be as in Proposition \ref{prop: seif 4-mfds}.
    Then $\Arf(K, Y, \s_0) = 0$, where $\s_0$ is the restriction of the unique spin structure on $B^5$ to $Y$.
\end{prop}
\begin{proof}
    The key observation is that the construction in Proposition \ref{prop: equi seif 4-mfds} gives a characteristic filling of $(K, Y, \s_0)$.
    Since $X$ embeds in $S^5$, $X$ is necessarily spin, and so since $[D] = 0$ in the relative homology of $X$, $D$ is characteristic in $X$.
    The fact that $X$ is spin also means that $\sigma(X) = 8\mu(Y, \s_0)$.
    In fact, spin structure on $Y$ induced by inclusion into $B^5$ extends over $X$ because $X$ is also embedded in $B^5$, and so the canonical spin characterization of $D$ in $X$ extends the canonical spin characterization of $K$ in $Y$ as they are both induced by the same ambient spin structure.
    Finally, since $K$ is null-homologous in $Y$, it has a unique Seifert framing in $Y$, and so the fact that $D$ is null-homologous rel. boundary implies that its normal Euler number relative to the Seifert framing is 0.
    We may therefore compute:
    \[ \Arf(K, Y, \s_0) = \frac{D \cdot D - \sigma(X)}{8} + \mu(Y, \s)  = \frac{0}{8} - \frac{\sigma(X)}{8} + \frac{\sigma(X)}{8} = 0.\]
\end{proof}




\section{Constructing and obstructing concordances} \label{sec: ann conc}

Following Proposition \ref{prop: seif 4-mfds give barf 0}, we can now define our concordance invariant for Montesinos annuli.

\subsection{The concordance invariant for annuli} \label{subsec: alpha}

\begin{defn}\label{def: alpha}
    Let $N$ be a Montesinos annulus in $S^1 \times S^3$, and let $S$ be its corresponding sphere in $S^4$ with removed sphere $R$.
    Let $Y$ be an Seifert 3-manifold for $S$ with $Y \cap R = A$ a null-homologous arc in $Y$.
    We define the invariant $\alpha$ of $N$ to be:
    \[ \alpha(N) := \Arf(A, Y, \s_0) \in \Z/2\Z,\]
    where $\s_0$ denotes the spin structure on $Y$ induced by the unique spin structure on $B^5$ via inclusion.
\end{defn}

The goal of this section is to prove the following classification theorem.

\begin{thm} \label{thm: annuli classification}
    There exists a $\Z/2\Z$-valued, additive relative concordance invariant $\alpha$ for Montesinos annuli $N$ in $S^1 \times B^3$ such that:
    \begin{enumerate}
        \item if $\alpha(N) = 0$, then $N$ is slice,
        \item if $\alpha(N) = 1$, then $N$ is concordant to the one-twist spun trefoil minus its spin axis, and
        \item every Montesinos annulus in $S^1 \times B^3$ satisfies exactly one of the above two conditions.
    \end{enumerate}
    Thus the smooth concordance group of Montesinos annuli $\calC _{2} ^\circ$ is isomorphic to $\Z/2\Z$.
\end{thm}

We start by showing that $\alpha$ is a well-defined quantity independent of choice of $Y$ and $A$, for which we need the following technical lemma.

\begin{lem} \label{lem: pushing into b5}
    Let $S$ be a 2-knot in $S^4$ and let $R$ be a smoothly unknotted 2-sphere in $S^4$ such that $S \cap R$ consists of two points.
    Let $Q$ be a smoothly unknotted 3-sphere embedded in $S^5$ such that $Q \cap S^4 = R$.
    Let $Y$ be a Seifert 3-manifold for $S$ such that $Y \cap R$ consists of a null-homologous arc $A$, and consider $Y$ as the image of an embedding $\psi \colon M^3 \to S^4$ for some 3-manifold with spherical boundary $M$. 
    Then $\psi$ is isotopic rel. boundary to a proper embedding $\psi'$ of $M$ in $B^5$ such that there exists a Morse function $f$ on $B^5$ with respect to which $\psi'(M)$ is in ordered Morse position.
    Furthermore, the arc $A' \colon = \psi'(M) \cap Q$ can be taken to be fully vertical with respect to $f|_{\psi'(M)}$ outside of a single critical point.
\end{lem}
\begin{proof}
    Choose a relative Morse function $g: M \to [0,1]$ such that:
    \begin{enumerate}
        \item $g|_{\partial M} \equiv 0$,
        \item $\psi(g^{-1}(1))$ is a single point in $A$, and
        \item away from $g^{-1}(1)$, $\psi^{-1}(A)$ is fully vertical with respect to $g$.
    \end{enumerate}
    Fix a local Morse function $f$ on a half-ball $U$ containing $Y$ in ${B^5}$ such that $f\vert_{\partial B^5} \equiv 0$ and $Q$ is fully vertical with respect to $f$.
    (Parametrizing $B^5$ and $Q$ as in Proposition \ref{prop: seif 4-mfds} will do.)
    We can then construct an isotopy $\Psi\colon M \times [0,1] \to f^{-1}([0,1])$ given by:
    \[ \Psi (x, t) = (\psi(x), t\cdot g\circ \psi(x)) \]
    
    Then $\Psi(x,0) = \psi(x)$ in the bounding ${S^4}$, and $\Psi(x,1)$ is a proper embedding $Y'$ of $M$ in ${B^5}$ such that $f|_{A'}$ for $A' := Y' \cap Q$ has exactly one maximum and is otherwise fully vertical with respect to $f|_{Y'}$ by our choice of $g$ and $\Psi$.
    Here again we use Proposition \ref{prop: omp} to rearrange $Y$ via ambient isotopy to ordered Morse position without changing the single critical point of $f|_{A'}$ or the full verticality of $A'$ with respect to $f$ away from this critical point.
\end{proof}

This lemma shows that it is possible to perturb Seifert 3-manifolds into the 5-ball in such a way that their intersections with the removed ball $Q$ is fully vertical with respect to the restriction of a Morse function $f$ to $Y$. 
This will enable us to use the constructions in Propositions \ref{prop: seif 4-mfds} and \ref{prop: seif 4-mfds give barf 0} to show that the $\alpha$ invariant is well-defined.

\begin{prop} \label{prop: barf is invt of N}
    The value of $\alpha(N)$ is independent of the choice of $Y$ and $A$.
\end{prop}
\begin{proof}
    Let $S$ be the corresponding sphere for $N$.
    Let $Y_1$ and $Y_2$ be two Seifert 3-manifolds for $S$, each intersecting the removed sphere $R$ in $S^4$ in null-homologous arcs $A_1$ and $A_2$.
    By the abuse of notation, in this proof we will speak of isotoping $Y_i$ and $A_i$ rel. boundary.
    $A_i$ always refers to the intersection of the 3-manifold $Y_i$ with $Q$, and is always topologically an arc.

    Consider the $S^4$ in which $S$ is embedded as the equator of a copy of $S^5$, such that this $S^4$ bounds two copies $B_1$ and $B_2$ of $B^5$ on either side.
    In this setup, we view $R$ as the equator of a smoothly unknotted 3-sphere $Q$ in $S^5$.

    Parametrize a 5-dimensional neighborhood $U$ of $S$ with coordinates $(x, y, z, w, t)$ so that the $S^4$ containing $S$ intersects $U$ in the subspace defined by setting $t = 0$.
    Then by Lemma \ref{lem: pushing into b5}, for $i = 1,2$ we may isotope $Y_i$ rel. boundary into $B_i$ such that each each arc of intersection $A_i$ of $Y_i$ with $Q$ is fully vertical with respect to the function $f \colon U \to \R$ given by projection onto the $t$-axis when restricted to $Y_i$.
    By Proposition \ref{prop: omp}, may then isotope each $Y_i$ rel. boundary to ordered Morse position with respect to $f$ such that this isotopy is supported away from $A_i$, maintaining the full verticality of the $A_i$'s with respect to $f|_{Y_i}$.

    Let $Y$ be the union of $Y_1$ and $Y_2$ along $S$, and $K$ be the union of $A_1$ and $A_2$ along their boundary.
    Then $Y$ is a closed 3-manifold in $S^5$ which intersects $Q$ in a knot $K$ which is fully vertical with respect to $f|_Y$ except for exactly two critical points and is null-homologous in $Y$, since each of the $A_i$'s was null-homologous in its $Y_i$.
    Thus $Y$, $A$, and $Q$ satisfy the conditions of Proposition \ref{prop: seif 4-mfds give barf 0}, and hence $\Arf(K,Y,\s_0) = 0$ where $\s_0$ is the restriction of the ambient spin structure on $S^5$ to $Y$.
    But $(K, Y, \s_0)$ is the connected sum of pairs $(A_1, Y_1,\s_0) \#(A_2, Y_2, \s_0)$, and so by Lemma \ref{lem: barf additivity}, we have that $\Arf(A_1, Y_1, \s_0) = \Arf(A_2, Y_2, \s_0)$.
\end{proof}

We similarly deduce the following proposition, which is one direction of the classification of Montesinos annuli up to concordance.
The idea is to use the same trick of building a 3-manifold intersecting a 3-sphere in a fully vertical knot as in the previous proposition, but this time, one of the two 3-manifolds glued together is a slice ball for $S$ which is not necessarily boundary parallel, and thus can't on its own be used to compute $\alpha$.

\begin{prop} \label{prop: slice means alpha 0}
    Let $N$ be an annulus in $S^1 \times B^3$ and suppose that $N$ is slice. 
    Then $\alpha(N) = 0$.
\end{prop}
\begin{proof}
    Parametrize a neighborhood of $N$ in $S^5$ with coordinates as in the previous proposition.
    Let $f$ be a local Morse function on this neighborhood which agrees with projection onto the $t$-coordinate.
    
    Since $N$ is slice, there exists a smooth embedding $\psi$ of $S^1 \times I \times I$ into $S^1 \times B^3 \times I$ such that $\psi(-,0) = N$, $\psi(-,1) = N_U$ the unknotted Montesinos annulus, and $\psi$ is level-preserving on $S^1 \times S^0 \times I \subset S^1 \times I \times I$.
    We may therefore embed $S^1 \times B^3 \times I$ into $S^5$ such that the final $I$ factor agrees with the $t$-coordinate in our parametrization.
    Under this embedding, the vertical boundary $S^1 \times S^0 \times I$ of the concordance is fully vertical with respect to $f$.
    Thus the union $C$ of the corresponding spheres for the annuli in each $t$-slice intersect the union of the removed spheres in each $t$-slice in two arcs $A_+$ and $A_-$ which are fully vertical with respect to $f|_C$.
    Let $Q$ be an unknotted 3-sphere in $S^5$ such that its intersection with our parametrized neighborhood is the $(z,w,t)$-subspace, and note that $Q \cap C$ is precisely $A_+ \cup A_-$.
    
    Now since $N_U$ is unknotted, we may find an unknotted 3-ball $B$ bounded by its corresponding sphere $S_U$ in the $t = 1$ slice which intersects $Q$ in an unknotted arc.
    By Lemma \ref{lem: pushing into b5}, we may find an isotopy of $B$ rel. boundary to some 3-ball $B' \subset S^5$ such that $B' \cap Q$ is an unknotted arc which is fully vertical with respect to $f|_{B'}$.

    Now let $D$ be the 3-ball in $S^5$ given by the union $B' \cup C$.
    Then $\partial D$ is the corresponding sphere $S$ for $N$ embedded in the $t=0$ slice.
    Note that since $A_+, A_-$, and $B' \cap Q$ were all fully vertical with respect to $f$, $D \cap Q$ is a single arc $\ell$ which is unknotted in $D$ and is fully vertical with respect to $f|_D$.

    Finally, let $Y$ be a Seifert 3-manifold for $S$ such that $Y$ intersects $Q$ in a null-homologous arc $A$, and isotope $Y$ into the $t \leq 0$ half of $S^5$ by Lemma \ref{lem: pushing into b5} such that its intersection with $Q$ is fully vertical with respect to $f|_Y$.
    The union $Y \cup D$ is then a closed 3-manifold embedded in $S^5$ such that its intersection $A \cup \ell$ with $Q$ is a null-homologous knot which is fully vertical with respect to $f|_{Y \cup D}$ apart from two critical points.
    We may then apply Proposition \ref{prop: seif 4-mfds give barf 0} to conclude $\Arf(A \cup \ell, Y \cup D, \s_0) = 0$, and hence by the additivity of the Arf invariant established in Lemma \ref{lem: barf additivity}, we have that $\Arf(A, Y, \s_0) = \Arf(\ell, D, \s_0)$.
    But since $D$ is a 3-ball and $\ell$ is unknotted, $\Arf(\ell, D, \s_0) = 0$, and so $\alpha(N) := \Arf(A, Y, \s_0) = 0$ as well.
\end{proof}

Finally, we establish the additivity of the $\alpha$ invariant, which follows quickly from the additivity of the Arf invariant.

\begin{lem} \label{lem: alpha additivity}
    Let $N$ and $N'$ be Montesinos annuli. Then $\alpha(N \natural N') = \alpha (N) + \alpha (N')$.
\end{lem}
\begin{proof}
    Let $Y$ and $Y'$ be Seifert 3-manifolds for the corresponding spheres $S$ and $S'$ and let $A$ and $A'$ be null-homologous arcs of intersection of $Y, Y'$ with removed spheres $R, R'$.
    Then by the stacking construction given in Definition \ref{def: stacking}, we see that the partial boundary connected sum $Y$ of $Y$ and $Y'$ forms a Seifert 3-manifold for $N \natural N'$ which intersects the connected sum of the removed spheres in the boundary connected sum of the tangles $A$ and $A'$.
    The result then follows from Lemma \ref{lem: barf additivity}.
\end{proof}

\subsection{Constructing slice balls} \label{subsec: slice balls}

The final step in the classification of Montesinos annuli up to concordance is to actually build the concordances.
We do this by building concordances of the corresponding spheres while continuously keeping track of their intersection with the removed 3-ball $Q$.
We will work in the same parametrized neighborhood as used in the constructions in Section \ref{sec: seif mfds}.

The general strategy of the construction is the same as that of Kervaire and Sunukjian: build a Seifert manifold $Y$ for $S$, push it into the interior of the 5-ball, and then ambiently surger it inside of $B^5$ to a 3-ball.
However, because we want our final slice ball $B$ for $S$ to intersect $Q$ in an unknotted arc $\ell$ so that we may remove $Q$ to obtain a slice cylinder for $N$, we must keep track of the intersection of $Y$ with $Q$ at each stage of the surgery.
This imposes an extra restriction that differs from the usual obstruction to ambient surgery in dimension 5, discussed in Subsubsection \ref{subsubsec: avoid Q}.

\begin{figure}[h]
    \labellist
    \small\hair 2pt
    \pinlabel {$Y$} at 180 160
    \pinlabel {$Q$} at 100 105
    \pinlabel {$A$} at 230 252
    \pinlabel {$S$} at 317 380
    \pinlabel {$R$} at 50 350
    \pinlabel {$S^4$} at 500 370
    \pinlabel {$B^5$} at 490 180 
    \endlabellist
    \centering
    \includegraphics[width=0.7\linewidth]{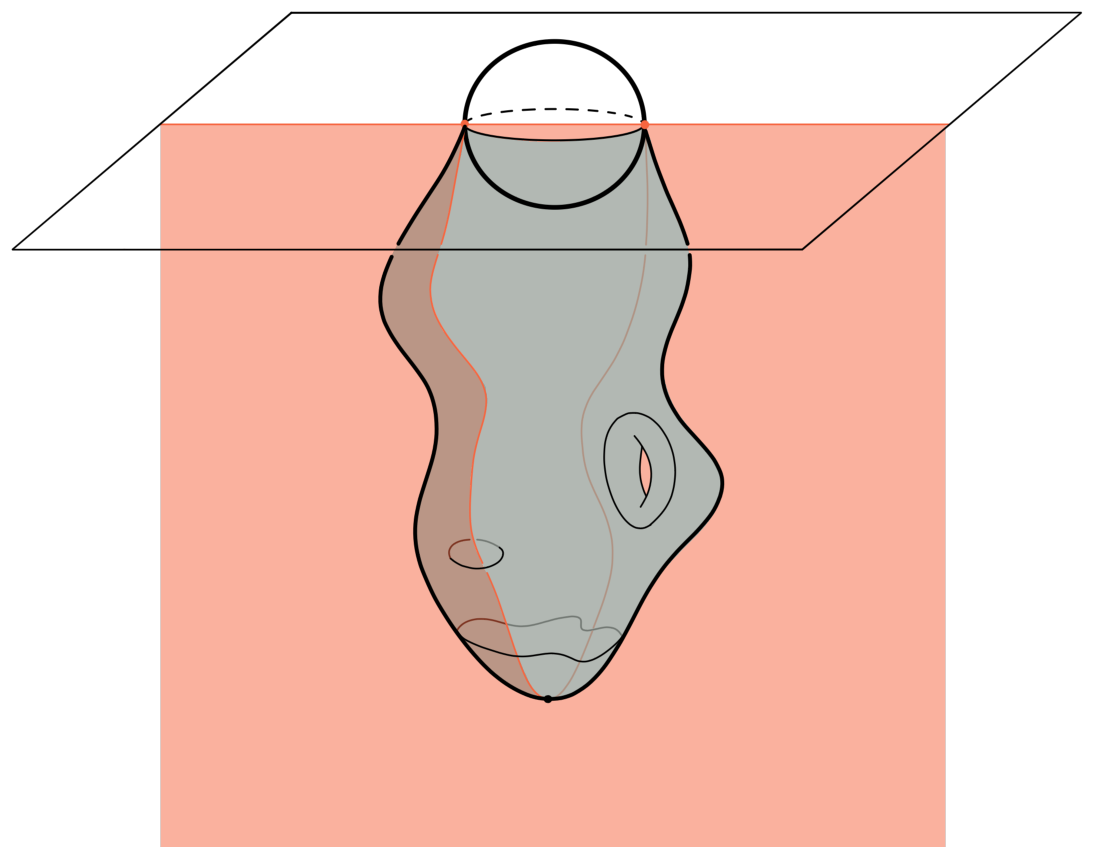}
    \caption{A schematic of the setup in preparation for surgery. Example surgery curves are shown.}
    \label{fig:placeholder}
\end{figure}

\subsubsection{Spin structures and ambient surgery in dimension 5} \label{subsubsec: spin surgery}

We first describe a restraint on ambient surgeries for 3-manifolds in dimension 5.
We adapt the language of Sunukjian \cite{sun15}, although this obstruction was known to Kervaire and overcome in his original proof of sliceness of 2-knots.

Recall that a \textit{trivialization} $\xi$ of a vector bundle $E$ over base $M^n$ with fiber $F^k$ is an identification of $E$ with $M \times F$.
One can specify a trivialization directly with a map, or with a choice of $k$ linearly independent sections of $E$.
We use the terms trivialization and \textit{framing} interchangeably.
Recall also that a spin structure on an $n$-manifold $M$ can be thought of as a trivialization of the tangent bundle $TM$ restricted to the 1-skeleton $M_1$ that extends over the 2-skeleton $M_2$.
A discussion of the equivalence of this definition to the usual one involving lifts of principal $\SO(n)$-bundles can be found in  \cite[Ch. VI]{kir89}.

Now suppose that $Y$ is a 3-manifold, and let $\mathfrak{s}$ be a spin structure on $Y$.
Let $\gamma$ be a curve in $Y$.
We say that $\xi$-framed surgery along $\gamma$ is a \textit{spin surgery} if the 4-manifold $X$ obtained by attaching a 4-dimensional 2-handle along $\gamma$ with framing $\xi$ to $Y \times \{1\}$ inside $Y \times [0,1]$ has a spin structure $\mathfrak{t}$ such that $\mathfrak{t}\vert_{Y \times \{0\}} = \mathfrak{s}$.

Finally, we recall the definition of ambient surgery. 
Let $Y$ be a 3-manifold embedded in $B^5$, and let $\gamma$ be a curve in $Y$ such that $\gamma$ bounds a disk $\Delta$ into $B^5 \setminus Y$.
Let $E$ be the total space of an embedded 2-disk sub-bundle of the normal bundle $\nu_{B^5} (\Delta)$ such that $\partial E \cap Y = \nu_Y (\gamma)$.
Lastly, let $\xi$ be a framing of $\nu_Y (\gamma)$ that extends to the unique trivialization of $E$ over $\Delta$.
Then the manifold $Y'$ embedded in $B^5$ given by $Y' := Y \setminus \nu_Y(\gamma) \cup (\Delta \times \partial D^2)$ is the result of \textit{ambient surgery} along $\gamma$ with framing $\xi$.

It is the purpose of the following lemma to determine when ambient surgery can be performed; that is, given $\gamma$ and $\Delta$, for which framings $\xi$ we can perform this operation.

\begin{lem} \label{lem: ambient iff spin}
    Let $Y$ be a 3-manifold embedded in $B^5$, and let $\gamma$ be a curve in $Y$ with a choice of framing $\xi$ such that $\gamma$ bounds a disk $\Delta$ into $B^5 \setminus Y$.
    Let $\mathfrak{s}_0$ be the unique spin structure on $B^5$.
    Let $Y'$ be the result of $\xi$-framed surgery along $\gamma$.
    Then the surgery $Y$ to $Y'$ can be realized as an ambient surgery in $B^5$ if and only if $\xi$-framed surgery along $\gamma$ is a spin surgery with respect to $\mathfrak{s}_0\vert_Y$.
\end{lem}

\begin{proof}
    Since $\gamma$ is a curve in an orientable a 3-manifold, trivializations of its tubular neighborhood $\nu_Y(\gamma)$ are in bijection with $\pi_1\SO(2) \cong \Z$ and carry an action of $\Z$ which increases or decreases the number of full counterclockwise twists.
    On the other hand, trivializations of $\nu_{B^5}(\gamma)$ are in bijection with $\pi_1\SO(4) \cong \Z/2\Z$.
    Indeed, the inclusion map of $\nu_Y (\gamma)$ into $\nu_{B^5}(\gamma)$ induces the mod 2 map from $\Z$ to $\Z/2\Z$.
    Thus, if $\xi$ and $\xi'$ are two different trivializations of $\nu_Y (\gamma)$ that differ by an even number of twists, then when included into $\nu_{B^5}(\gamma)$, they become homotopic.
    
    Since $\Delta$ is contractible, every 2-disk bundle over it is canonically trivial; denote this trivialization $\zeta$.
    Thus given a 2-disk bundle over $\Delta$, we expect that half of the possible framings on $\gamma$ can be extended over the bundle, depending on their parity relative to $\zeta$.
    
    Now $\Delta$ is embedded in $B^5$, and so the unique trivialization of its normal bundle $\nu_{B^5}(D)$ is exactly the unique spin structure $\mathfrak{s_0}$ on $B^5$ restricted to $\Delta$.
    This then further restricts to a trivialization $\zeta$ on any 2-disk sub-bundle of $\nu_{B^5}(\Delta)$, and in particular on $\nu_{B^5}(\partial \Delta)$, which is just $\nu_{B^5}(\gamma)$.
    Thus the framings $\xi$ on $\nu_{Y}(\gamma)$ for which there exist embedded $D^2$ sub-bundles of $\nu_{B^5}(\Delta)$ over which they can be extended are precisely those that agree with $\zeta$ restricted to $\nu_{B^5}(\gamma)$.
    Since $\zeta$ was induced by the restriction of the ambient spin structure on $B^5$, we see that surgery along $\gamma$ with framing $\xi$ must be spin surgery. 
\end{proof}

\subsubsection{Avoiding the removed 3-ball} \label{subsubsec: avoid Q}

Here we address the restraint on the allowable surgeries imposed by requiring that a slice ball we construct for $S$ actually give a slice cylinder for the corresponding annulus $N$.

\begin{figure}[h]
    \labellist
    \small\hair 2pt
    \pinlabel {$Y$} at 180 160
    \pinlabel {$Q$} at 100 105
    \pinlabel {$A$} at 230 252
    \pinlabel {$S$} at 317 380
    \pinlabel {$R$} at 50 350
    \pinlabel {$S^4$} at 500 370
    \pinlabel {$B^5$} at 490 180 
    \endlabellist
    \centering
    \includegraphics[width=0.7\linewidth]{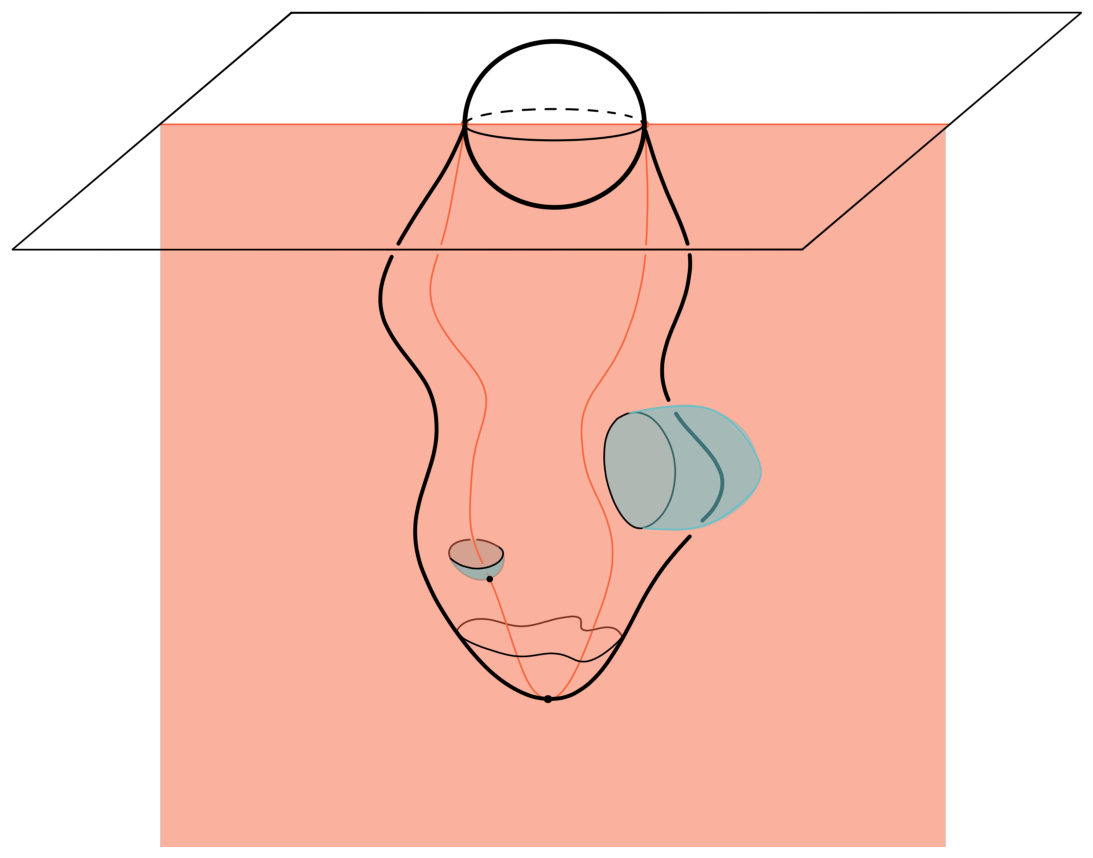}
    \caption{A schematic of two disks bounded by curves in $Y$. Note that the disk on the right is disjoint from $Q$, but the disk on the left intersects it transversely in a single point.}
    \label{fig: y with surg disks}
\end{figure}

\begin{figure}[h]
    \labellist
    \small\hair 2pt
    \pinlabel {$\gamma$} at -10 100
    \pinlabel {$A$} at 23 40
    \pinlabel {$Q$} at 230 20
    \pinlabel {$\gamma$} at 230 120
    \pinlabel {$H$} at 312 60
    \endlabellist
    \centering
    \includegraphics[width=0.5\linewidth]{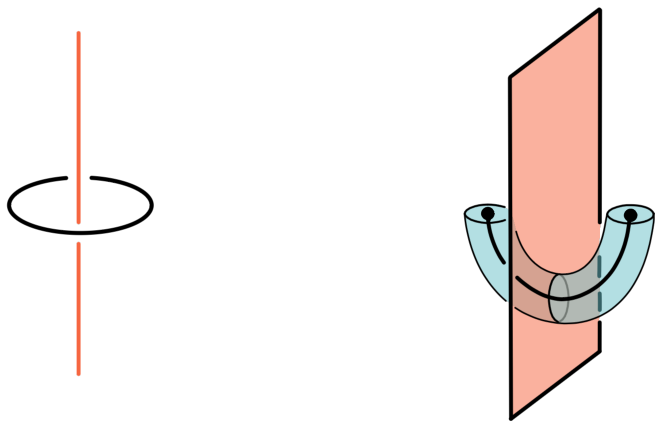}
    \caption{New components of $Y \cap Q$ introduced by errant handle attachment. 
    On the left, we see an example surgery curve $\gamma$ linking $Q$ inside $Y$ (and thus linking their intersection, which is $A$). 
    On the right, we see a different 3-dimensional projection of the same intersection, in which $\gamma$ appears as a pair of points and $Q$ as a plane. 
    The arc connecting these two points and intersecting $Q$ is one slice of the disk $\Delta$ that we want to surger along, with the 2-disk bundle $H$ shaded. 
    However, as we see, this neighborhood intersects $Q$ in a disk, and so surgery along $\gamma$ will introduce a new component of $Y \cap Q$ given by the intersection of the boundary of $H$ with $Q$.}
    \label{fig: new intersect}
\end{figure}

Suppose we wished to surger along a curve $\gamma$ in $Y$ bounding a disk $\Delta$ in $B^5$ such that $\Delta \cap Q$ was nonempty.
Generically, $\Delta$ intersects $Q$ in points, so let's examine one such intersection point $q$.
When $Y$ is replaced by $Y' = (Y \setminus \nu(\gamma)) \cup_{\partial(\nu(\gamma))} \partial(H)$ for $H$ some 2-disk sub-bundle of $\nu(\Delta)$ thought of as a 4-dimensional 2-handle, $Y' \cap Q$ now includes the belt sphere of $H$ over $q$.
This is illustrated in Figures \ref{fig: y with surg disks} and \ref{fig: new intersect}.

This is a major problem for our surgery program, because once a new closed component of this intersection is introduced, it is not clear that it can be removed by ambient surgery.
We therefore want to restrict ourselves to surgery along curves in $Y$ which bound disks with interiors in $B^5 \setminus (Y \cup Q)$.
Thankfully, we can compute the fundamental group of this complement.

\begin{prop} \label{prop: pi1 z2}
    $\pi_1(B^5\setminus (Y \cup Q)) \cong \Z^2$.
\end{prop}
\begin{proof}
    We make a rising waters argument; compare with \cite{gs00}, Prop. 6.2.1 for this type of argument for surfaces in 4-manifolds.
    Consider the radial Morse function $f$ on $B^5$ from Lemma \ref{lem: pushing into b5}.
    Recall that both $Y$ and $Q$ are properly embedded and boundary parallel, and their intersection is a single properly embedded arc $A$.
    The key observation is that by construction, $Y, Q,$ and $A$ each have a single 0-handle in the handle decomposition induced by the restriction of $f$.
    
    The 0-handles of $Y$ and $Q$ each contribute a 1-handle each to the complement.
    The minimum point of their intersection then contributes a 2-handle, the attaching sphere of which winds around the commutator of the meridians of both $Y$ and $Q$ near the intersection point.
    Thus the fundamental group of the sublevel set of $f$ just above the first intersection point is $\Z^2$.
    However, all subsequent critical points of $Y \cup Q$ are in $Y$ and are of index 1 or 2, which contribute 2- or 3-handles to the complement.
    Of these, only the 2-handles can contribute to the fundamental group relations of $B^5 \setminus (Y \cup Q)$, but these only introduce extra commutation relations between meridians that already commute.
\end{proof}

Hence we conclude that if $\gamma$ is a curve in $Y$ which is null-homologous in $B^5 \setminus Q$, then it bounds a disk in $B^5 \setminus (Y \cup Q)$ and we may therefore surger along it (with a spin framing) without introducing new components of $Y \cap Q$.
Note that this condition is satisfied it $\gamma$ is null-homologous in $Y \setminus A$.

\subsubsection{Finding surgery curves}

The discussions of the previous two subsubsections amount to the following rule when performing ambient surgery in this setup: 
we may surger $Y$ along any curve $\gamma$ with framing $\xi$ so long as:
\begin{enumerate}
    \item $\xi$-framed surgery is spin surgery with respect to the restriction of the ambient spin structure $\s_0$ on $B^5$ to $Y$, and 
    \item $\gamma$ must be null-homologous in $Y \setminus A$.
\end{enumerate}
For ease of exposition, we will call pairs $(\gamma, \xi)$ satisfying these properties \textit{viable}, and surgeries of $Y$ along them \textit{viable surgeries}.

Our invariant $\alpha$ was chosen precisely so as to be invariant under viable surgeries:
\begin{prop} \label{prop: alpha under surgeries}
    Let $Y$ be a Seifert manifold for $S$ intersecting $Q$ in a null-homologous arc $A$, and let $Y'$ be a 3-manifold properly embedded in $B^5$ bounded by $S$ intersecting $Q$ in null-homologous arc $A'$.
    Suppose that $Y'$ is related to $Y$ by viable surgeries.
    Then $\Arf(A, Y, \s_0) = \Arf(A', Y', \s_0)$.
\end{prop}
\begin{proof}
    Recall that abstractly, we may cap off each pair of arcs and 3-manifolds with trivial arcs in a 3-ball, and that we use these closed versions to compute the Arf invariant, as in Definition \ref{def: barf arc}.
    For ease of argument, continue to call the capped-off manifolds $Y$ and $Y'$, and the knots inside them $A$ and $A'$.
    
    A collection of viable surgeries $(\gamma_i, \xi_i)$ relating $Y$ and $Y'$ gives a recipe for a smooth characteristic concordance in a spin 4-manifold $X$ (in the sense of Definition \ref{def: char conc}) between $(A, Y, \s_0)$ and $(A', Y', \s_0)$. 
    Specifically, because viable surgery curves are required to be null-homologous in $Y\setminus A$, the product $A \times I$ embeds in $X$ as a characteristic annulus cobounded by $A$ and $-A'$.
    Thus by Proposition \ref{prop: char conc iff arf}, $\Arf(A, Y, \s_0) = \Arf(A', Y', \s_0)$.
\end{proof}

With these restrictions in mind, we turn to finding viable surgery curves and disks, which we do in two steps.
The first is to ignore $A$ and simply find surgery curves that realize an embedded spin cobordism between $Y$ and $B^3$ inside $B^5$.
The second is to then modify the 3-ball we surgered $Y$ to in order to unknot its intersection with $Q$.
This is really the same procedure as the construction of spin characteristic fillings given in Proposition \ref{prop: spin fill}, but performed ambiently inside of $B^5$.

\begin{lem} \label{lem: Y to B3}
    There exists a sequence of viable surgeries from $Y$ to $B^3$ such that the resulting $B^3$ intersects $Q$ in a single arc.
\end{lem}
\begin{proof}
    Recall that by Kaplan's theorem \cite{kap79}, for any spin 3-manifold $(Y, \s)$, there exists a sequence of spin surgeries taking $(Y, \s)$ to the 3-ball with its unique spin structure $(B, \s_B)$.
    We use this to find starting curves and framings for our surgeries on $(Y, \s_0)$, since by Lemma \ref{lem: ambient iff spin} we can realize any spin surgery as an ambient surgery in $B^5$.
    Once we identify the curves we wish to surger $Y$ along, we isotope them inside $Y$ so that they are null-homologous in $B^5 \setminus (Y \cup Q)$, allowing them to pass through $A$ if necessary.
    Such an isotopy inside $Y$ can always be found as $A$ is null-homologous in $Y$ rel. boundary, and so if a surgery curve $\gamma$ in $Y$ does not bound a disk in $B^5 \setminus Y \cup Q$, we may use a relative Seifert surface for $A$ in $Y$ bounded by $A$ to isotope it along a path to $A$ and push it past $A$ to remove the intersection. 
    This is similar to the trick shown in Figure \ref{fig: 1d whitney}.
    The result of the surgery on these curves will be to turn $Y$ into a 3-ball $B$, realizing $A$ as a knotted arc in $B$.
\end{proof}

It remains now to unknot $A$ through viable ambient surgeries.
These surgeries must preserve the homeomorphism type of the 3-ball in which $A$ lives, and thus amount to changing the isotopy class of $A$ as an arc in $B^3$ by changing which $B^3$ in $B^5$ we think of $A$ as living in.

\begin{lem} \label{lem: unknot A}
    Suppose that $\Arf(A, B, \s_0) = 0$. Then there exists a sequence of viable surgeries on $B$ which realizes $A$ as an unknotted arc in a 3-ball.
\end{lem}
\begin{proof}
    Note that since $B$ is a 3-ball and therefore has a unique spin structure, $\Arf(A, B, \s_0)$ is equal to the classical Arf invariant of the knot $K$ given by capping off $B$ with another copy of $B^3$ and $A$ with a trivial arc.
    Recall by Lemma \ref{lem: all knots char slice} that if $\Arf(K) = 0$, then $K$ is band pass equivalent to the unknot, and that these band passes can be achieved by 0-framed (thus spin) surgery on curves which are null-homologous in the complement of $K$ in $S^3$, as shown in Figure \ref{fig: band pass}.
    These curves and coefficients then represent viable surgeries on $B$.
    We may take these curves to be contained in $B$, and since they are null-homologous in $B \setminus A$, by Proposition \ref{prop: pi1 z2}, they bound disks in the complement of $B \cup Q$ in $B^5$, and hence surgering along them does not introduce any new components of $B \cap Q$.
    The result of these surgeries is then a (possibly different) 3-ball bounded by $S$ which intersects $Q$ in an unknotted arc $\ell$.
\end{proof}

We can now put together the entire ambient surgery argument to construct slice cylinders.

\begin{prop} \label{prop: alpha 0 means slice}
    If $\alpha(N) = 0$, then $N$ is slice.
\end{prop}
\begin{proof}
    Let $N$ be a Montesinos annulus, let $S$ be its corresponding sphere, and let $R$ be its removed sphere.
    By Proposition \ref{prop: seif 3-mfds} and Lemma \ref{lem: pushing into b5}, $S$ bounds a Seifert 3-manifold $Y$ which can be pushed into $B^5$ such that its intersection with $Q$ is an arc $A$.
    By Lemma \ref{lem: Y to B3}, $Y$ can then be ambiently surgered to a copy $B$ of $B^3$ through a sequence of viable surgeries, thus preserving $Y \cap Q$.

    By the definition of the invariant and Proposition \ref{prop: alpha under surgeries}, the assumption that $\alpha(N) = 0$ implies that $\Arf(A, Y, \s_0) = 0$.
    Therefore by Lemma \ref{lem: unknot A}, there exists a sequence of viable surgeries on $B$ which turn $B \cap Q$ into an unknotted arc $\ell$ in a 3-ball $B'$.
    Since $B'$ is a 3-ball bounded by $S$ with $B' \cap Q = \ell$ an unknotted arc, by removing a neighborhood of $Q$, we obtain a slice cylinder for $N$.
\end{proof}

We may now show that the $\alpha$ invariant completely classifies Montesinos annuli up to concordance.

\begin{proof}[Proof of Theorem \ref{thm: annuli classification}]
    Together, Propositions \ref{prop: alpha 0 means slice} and \ref{prop: slice means alpha 0} show that $N$ is slice if and only if $N$ has vanishing $\alpha$ invariant.
    We can then deduce that all annuli with nonvanishing $\alpha$ invariant are concordant by exploiting the group structure of annuli under stacking. 
    If $N$ and $N'$ are two Montesinos annuli such that $\alpha(N) = \alpha(N') = 1$, then $\alpha(N \natural N') = 0$ by Lemma \ref{lem: alpha additivity}, and so by the above construction, $N \natural N'$ is slice.
    We therefore have that every annulus with nonvanishing $\alpha$ invariant is the inverse of every other annulus with nonvanishing $\alpha$ invariant, and so since these annuli cannot represent the identity element in $\calC_{2} ^\circ$ by Proposition \ref{prop: slice means alpha 0}, they must all belong to the same concordance class which has order 2 in the concordance group.
    
    To see that $N_T$ generates $\calC_{2} ^\circ$, recall that in Example \ref{ex: N_T} we constructed a 3-ball bounded by $S^T$ which intersects the removed sphere in a trefoil tangle.
    Since the trefoil has Arf invariant 1, we have that $\alpha(N_T) = 1$, and so it cannot be slice by Proposition \ref{prop: slice means alpha 0}.
\end{proof}


\section{Concordance of periodic 2-knots} \label{sec: equi conc}

To port over the classification given in Theorem \ref{thm: annuli classification} to the setting of periodic 2-knots, we take a moment to recall the relationship between periodic 2-knots and Montesinos annuli.
Let $\tilde S$ be a $d$-periodic 2-knot in $S^4$.
Then we may take the quotient of $\tilde S$ by $\rho$ to obtain a new sphere $S$ inside $S^4$ whose intersection with $\Branch(\rho)$ consists of exactly two points.
Since $\Branch(\rho)$ is an unknotted 2-sphere, we may remove a neighborhood of it to obtain a Montesinos annulus $N$ in $S^1 \times B^3$ given by $S \setminus \nu(\{p_+, p_-\})$.
We call the annulus $N$ obtained in this way the \textit{quotient annulus} of $\tilde S$.
The result we would like to prove is:

\begin{thm} \label{thm: equi conc}
    For each positive integer $d$, there exists a $\Z/2\Z$-valued, additive concordance invariant $\alpha$ for $d$-periodic 2-knots $\tilde S$ such that:
    \begin{enumerate}
        \item if $\alpha_d(\tilde S) = 0$, then $\tilde S$ is $\Z/d\Z$-equivariantly slice,
        \item if $\alpha_d(\tilde S) = 1$, then $\tilde S$ is $\Z/d\Z$-equivariantly concordant to the $d$-twist spun trefoil, and
        \item every $d$-periodic 2-knot satisfies exactly one of the above two conditions.
    \end{enumerate}
    Thus the smooth equivariant concordance group $\calC_2 ^d$ is isomorphic to $\Z/2\Z$ for all $d$.
\end{thm}

To prove this, we will need to successively reprove many results from earlier in this paper in the equivariant setting by passing from equivariant objects to their quotients, applying earlier results, and then passing back by taking an appropriate cyclic branched cover.

\subsection{Equivariant Seifert manifolds}

The first step of this process is to construct equivariant Seifert manifolds for simply periodic 2-knots and 3-manifolds.
The relationship between these results and those in Section \ref{sec: seif mfds} is simply that the ambient subspaces $R$ and $Q$ mentioned in those proofs here appear as the images of the fixed sets of $\rho$ under the quotient map.
Recall our convention that objects which carry a cyclic group action are denoted with tildes, and objects which are invariant under the action are denoted with hats.
Objects in the quotient have no decoration.

\begin{prop} \label{prop: equi seif 3-mfds}
    Let $\tilde S$ be a periodic 2-knot in $\tilde S^4$.
    Then $\tilde S$ bounds a $\rho$-equivariant Seifert 3-manifold $\tilde Y \subset \tilde S^4$ such that $\Fix(\rho) \cap \tilde Y$ consists of an arc $A$ that is null-homologous in $\tilde Y$ rel. boundary.
\end{prop}
\begin{proof}
    Let $S$ be the quotient sphere for $\tilde S$ inside the copy of $S^4$ obtained after taking the quotient by $\rho$.
    Then $\Branch(\rho)$ is an unknotted 2-sphere intersecting $S$ in exactly two points.
    We may therefore apply Proposition \ref{prop: seif 3-mfds} to obtain a Seifert 3-manifold $Y$ for $S$ that intersects $\Branch(\rho)$ in an arc $A$ which is null-homologous in $Y$ rel. boundary.
    Since $A$ is contained in $\Branch(\rho)$, after taking the $d$-fold cyclic branched cover of $S^4$ along $\Branch(\rho)$, $Y$ lifts to an equivariant Seifert 3-manifold $\tilde Y$ whose intersection with $\Fix(\rho)$ is the lift of $A$.
    Finally, since $A$ was null-homologous in $Y$ rel. boundary, its lift is null-homologous in $\tilde Y$ rel. boundary.
    \end{proof}

The construction of equivariant Seifert 4-manifolds below is much the same.

\begin{prop} \label{prop: equi seif 4-mfds}
    Let $\tilde Y$ be a closed, simply periodic 3-manifold smoothly embedded in $\tilde S^5$, and suppose that $\tilde Y \cap \Fix(\rho)$ is a null-homologous knot $\hat K$ in $\tilde Y$.
    Suppose that there exists an equivariant Morse function $\tilde f$ on $\tilde S^5$ such that $\tilde Y$ is in ordered equivariant Morse position with respect to $\tilde f$, $\tilde f|_{\hat K}$ has exactly two critical points, and away from these two critical points, $\hat K$ is fully vertical with respect to $\tilde f|_{\tilde Y}$.
    Then $\tilde Y$ bounds a $\rho$-equivariant Seifert 4-manifold $\tilde X$ smoothly embedded in $\tilde S^5$ such that $\tilde X \cap \Fix(\rho)$ is a smoothly embedded disk $\hat D$ bounded by $\hat K$ which is null-homologous in $\tilde X$ rel. boundary. 
\end{prop}
\begin{proof}
    Let $Y$ be the quotient of $\tilde Y$ under $\rho$ inside $S^5$ such that $Y \cap \Branch(\rho)$ is the knot $K$ given by the image of $\hat K$ under the covering map.
    Since $\tilde f$ is equivariant, it also descends to a Morse function $f \colon S^5 \to \R$.
    Since $\tilde Y$ was in ordered equivariant Morse position with respect to $\tilde f$, $Y$ is in ordered Morse position with respect to $f$, and since $\hat K$ was contained in $\Fix(\rho)$, $f|_K$ has exactly two critical points and $K$ is otherwise fully vertical with respect to $f$.
    We therefore can apply Proposition \ref{prop: seif 4-mfds} to obtain a Seifert 4-manifold $X$ for $Y$ smoothly embedded in $S^5$ such that its intersection with $\Branch(\rho)$ is a smoothly embedded disk $D$ bounded by $K$ which is null-homologous in $X$ rel. boundary.
    Taking the branched cover of the whole setup over $\Branch(\rho)$ then recovers $\tilde Y$ with $\tilde Y \cap \Fix(\rho) = \hat K$ and the lifts of $X$ and $D$ give $\tilde X$ and $\hat D$ with the desired properties.
\end{proof}

\subsection{The equivariant concordance invariant} \label{subsec: alpha_d}

This subsection reproves the contents of Subsection \ref{subsec: alpha} in the context of periodic 2-knots.
For the following, fix a period $d$ for $\rho$.

\begin{defn}\label{def: alpha_d}
    Let $\tilde S$ be a $d$-periodic 2-knot in $\tilde S^4$.
    Let $\tilde Y$ be an equivariant Seifert 3-manifold for $\tilde S$ with $\tilde Y \cap \Fix(\rho) = \hat A$ a null-homologous arc in $\tilde Y$.
    Consider the manifold $Y$ and contained arc $A$ given by taking the quotient by $\rho$.
    We define the invariant $\alpha_d$ of $\tilde S$ to be:
    \[ \alpha_d(\tilde S) := \Arf(A, Y, \s_0) \in \Z/2\Z,\]
    where $\s_0$ denotes the spin structure on $Y$ induced by the unique spin structure on $B^5$ via inclusion.
\end{defn}

\begin{prop} \label{prop: barf is invt of S tilde}
    The value of $\alpha_d(\tilde S)$ is independent of the choice of $\tilde Y$ and $\hat A$.
\end{prop}
\begin{proof}
    Let $\tilde Y_1$ and $\tilde Y_2$ be two equivariant Seifert manifolds for $\tilde S$, and suppose that each intersects $\Fix(\rho)$ in arcs $\hat A_1$ and $\hat A_2$, respectively.
    Let $Y_1$ and $Y_2$ be their quotients by $\rho$, intersecting $\Branch(\rho)$ in arcs $A_1$ and $A_2$, respectively.
    As in the proof of Proposition \ref{prop: barf is invt of N}, parametrize a neighborhood $U$ of $S$ in $S^5$ with coordinates $(x,y,z,w,t)$ such that $\Branch(\rho) \cap U$ consists of the $(z,w,t)$-subspace, and that $S^4 \cap U$ is the $t=0$ subspace.
    We can then apply Lemma \ref{lem: pushing into b5} to isotope $Y_1$ and $Y_2$ rel. boundary into opposite copies $B_1$ and $B_2$ of $B^5$ bounded by $S^4$ such that each intersects $\Branch(\rho)$ in an arc, which we continue to denote $A_1$ and $A_2$.
    The lemma also provides us with a Morse function $f$, which we can take to be projection onto the $t$ coordinate, such that each $A_i$ is fully vertical with respect to $f|_{Y_i}$ away from a single critical point.
    By Proposition \ref{prop: omp}, we can isotope each $Y_i$ rel. boundary to ordered Morse position with respect to $f$ through isotopies supported away from the $A_i$'s, ensuring that they will lift to isotopies in the branched cover as in the proof of Proposition \ref{prop: oemp}.

    Let $Y$ be the union of $Y_1$ and $Y_2$ along $S$, and let $K$ be the union of $A_1$ and $A_2$ along their boundary.
    Then $Y, A,$ and $\Branch(\rho)$ satisfy the conditions of Proposition \ref{prop: seif 4-mfds give barf 0}, and so we have that $\Arf(K, Y, \s_0) = 0$ for $\s_0$ the unique spin structure on $S^5$.
    Conclude by the additivity of the Arf invariant in Lemma \ref{lem: barf additivity} that $\Arf(A_1, Y_1, \s_0) = \Arf(A_2, Y_2, \s_0)$.
\end{proof}

As in Proposition \ref{prop: slice means alpha 0}, we use the same trick to show that $\alpha_d$ vanishes on $\Z/d\Z$-equivariantly slice periodic 2-knots.

\begin{prop} \label{prop: equi slice means alpha_d 0}
    Let $\tilde S$ be a $d$-periodic 2-knot in $\tilde S^4$ and suppose that $\tilde S$ is $\Z/d\Z$-equivariantly slice. 
    Then $\alpha_d(\tilde S) = 0$.
\end{prop}
\begin{proof}
    Parametrize a neighborhood of the quotient sphere $S$ in $S^5$ with coordinates $(x,y,z,w,t)$ such that $S$ lies in the $t=0$ hypersurface, and such that $\Branch(\rho)$ consists of the $(z,w,t)$-subspace.
    Let $f$ be projection onto the $t$-axis.

    Since $\tilde S$ is equivariantly slice, there exists an equivariant concordance $\tilde C \subset \tilde S^4 \times I$ cobounded by $\tilde S$ and the unknot $\tilde U$.
    Recall that our definition of equivariantly slice requires this embedding to be level-preserving on the fixed sets of these actions.
    We may therefore embed the quotient concordance inside $S^4 \times I$ into our parametrized neighborhood of $S^5$ in such a way that $C\cap \Fix(\rho)$ is fully vertical with respect to $f|_C$.
    We may then cap off this quotient concordance with a spanning 3-ball $B$ bounded by the quotient unknot $U$ in the $t=1$ slice, pushed slightly deeper in the direction of increasing $t$ using Lemma \ref{lem: pushing into b5} such that its intersection with $\Branch(\rho)$ is fully vertical with respect to $f|_B$, to obtain a slice ball $B'$ for $S$ which intersects $\Branch(\rho)$ in an unknotted arc $\ell$ which is fully vertical with respect to $f|_{B'}$.

    Now let $\tilde Y$ be an equivariant Seifert 3-manifold for $\tilde S$ intersecting $\Branch(\rho)$ in a null-homologous arc, and let $Y$ be its quotient under $\rho$.
    Again by Lemma \ref{lem: pushing into b5}, we may isotope $Y$ into the $t \leq 0$ 5-ball such that its intersection with $\Branch(\rho)$ in an arc $A$ which is fully vertical with respect to $f|_Y$.
    Then the closed 3-manifold $Y \cup B'$ together with $A \cup \ell$ satisfies the conditions of Proposition \ref{prop: seif 4-mfds give barf 0}, and so we have that $\Arf(A \cup \ell, Y \cup B', \s_0) = 0$.
    Again by the additivity established in Lemma \ref{lem: barf additivity}, this implies $\Arf(A, Y, \s_0) = \Arf(\ell, B', \s_0)$, which then vanishes because $\ell$ is unknotted and $B'$ is a 3-ball.
\end{proof}

We can also establish additivity for $\alpha_d$.

\begin{lem} \label{lem: alpha_d additivity}
    Let $\tilde S$ and $\tilde S'$ be two $d$-periodic 2-knots inside two copies of the 4-sphere which equipped with $\Z/d\Z$ actions $\rho$ and $\rho'$, respectively. Then 
    \[\alpha_d(\tilde S \#_{\Z/d\Z} \tilde S') = \alpha_d (\tilde S) + \alpha_d (\tilde S').\]
\end{lem}
\begin{proof}
    Let $S$ and $S'$ be the quotient spheres of $\tilde S$ and $\tilde S'$.
    Let $Y$ and $Y'$ be Seifert 3-manifolds for the corresponding spheres $S$ and $S'$ and let $A$ and $A'$ be null-homologous arcs of intersection of $Y, Y'$ with $\Branch(\rho)$ and $\Branch(\rho')$, respectively.
    Then by the equivariant connected sum given in Definition \ref{def: equi conn sum}, we see that the partial boundary connected sum $Y$ of $Y$ and $Y'$ forms a Seifert 3-manifold for the quotient sphere of $\tilde S \#_{\Z/d\Z} \tilde S'$  which intersects the connected sum $\Branch(\rho) \# \Branch(\rho')$ in the boundary connected sum of the tangles $A$ and $A'$.
    The result then follows from Lemma \ref{lem: barf additivity}.
\end{proof}

\subsection{Building equivariant concordances}

We can now leverage the contents of Subsection \ref{subsec: slice balls} to construct equivariant slice balls for $d$-periodic 2-knots with vanishing $\alpha_d$ invariant.
The strategy to constructing the slice ball is as follows:
First, we obtain an equivariant Seifert 3-manifold $\tilde Y$ for $\tilde S$, and then take its quotient by $\rho$ to obtain some other 3-manifold $Y$ which is a Seifert manifold for $S$.
This realizes $\tilde Y$ as the cyclic cover of $Y$ branched over an arc $A$, where $\hat A:= \tilde Y \cap \Fix(\rho)$ descends to $A$ under the branched covering map.
Then just as in the non-equivariant case, we isotope $Y$ into the interior of ${B^5}$ and ambiently surger $Y$ to make it a 3-ball $B$ and then $A$ to an unknotted arc $\ell$.
In order to ensure that these surgeries can be performed ambiently without introducing new intersections with $\Branch(\rho)$ (and therefore lift to surgeries on the cover), we will require $\alpha_d$ to vanish.
We then lift the isotopy and surgery curves and disks to $\tilde Y \subset \tilde B^5$, obtaining a sequence of suitable equivariant surgeries taking $(\hat A, \tilde Y)$ to $(\hat \ell, \tilde B^3)$, since the latter is its own cyclic branched cover.

\begin{prop} \label{prop: alpha d zero means slice}
    Let $\tilde S$ be a $d$-periodic knot with $\alpha_d(\tilde S) = 0$.
    Then $\tilde S$ is $\Z/d\Z$-equivariantly slice.
\end{prop}
\begin{proof}
     Let $\tilde Y$ be an equivariant Seifert 3-manifold for $\tilde S$ intersecting $\Fix(\rho)$ in a null-homologous arc $\hat A$.
     Let $S$ be the quotient sphere for $\tilde S$, and note that $\tilde Y$ descends to a Seifert 3-manifold $Y$ for $S$ intersecting $\Branch(\rho)$ in a null-homologous arc $A$.
     By Lemma \ref{lem: pushing into b5}, we may isotope $Y$ rel. boundary into the interior of $B^5$ while keeping its intersection with $\Branch(\rho)$ still a single arc, which we continue to call $A$.
     We are now in precisely the same setup as in Subsection \ref{subsec: slice balls}, taking $\Branch(\rho)$ to be the removed sphere.
     We may therefore apply Lemma \ref{lem: Y to B3} to find a sequence of ambient surgeries taking $Y$ to a 3-ball $B$ bounded by $S$ which do not introduce any new intersections with $\Branch(\rho)$.

     As in Proposition \ref{prop: alpha 0 means slice}, the condition that $\alpha_d(S) = 0$ requires that $\Arf(A, Y, \s_0) = 0$, which by Proposition \ref{prop: alpha under surgeries} means that $\Arf(A, B, \s_0) =0$.
     We can then apply Lemma \ref{lem: unknot A} to find another sequence of ambient surgeries on $B$ turning it into another 3-ball $B'$ such that $B' \cap \Branch(\rho)$ is an unknotted arc $\ell$.
     Since the disks surgered along did not intersect the branch set of $\rho$, we may lift the surgery curves to a set of equivariant curves on $\tilde Y$.
     The framed disks bounded by these curves lift as well, and so we obtain from this procedure a set of equivariant ambient surgeries converting $\tilde Y$ into an equivariant 3-ball bounded by $\tilde S$.
\end{proof}

\begin{proof}[Proof of Theorem \ref{thm: equi conc}]
    The proof is exactly the same as in Theorem \ref{thm: annuli classification}.
    Propositions \ref{prop: equi slice means alpha_d 0} and \ref{prop: alpha 0 means slice} show that a $d$-periodic 2-knot $\tilde S$ is $\Z/d\Z$-equivariantly slice if and only if $\alpha_d(\tilde S) = 0.$
    Since $\alpha_d$ is additive under equivariant connected sum by Proposition \ref{lem: alpha_d additivity} and valued in $\Z/2\Z$, every periodic 2-knot which is not $\Z/d\Z$-equivariantly slice is the inverse of every other such 2-knot in the corresponding equivariant concordance group $\calC_2^d$. 
    Hence all such 2-knots must be equal and have order 2 in $\calC_2 ^d$, which thus must be isomorphic to $\Z/2\Z$.
\end{proof}

\bibliographystyle{alpha}
\bibliography{bib}

\end{document}